\documentclass[a4paper,reqno]{amsart}
\author{Julia Brandes}
\title[Forms representing forms and linear spaces on hypersurfaces]{Forms representing forms \\ and linear spaces on hypersurfaces}
\address{School of Mathematics, University of Bristol, University Walk, Clifton, Bristol BS8 1TW, United Kingdom}
\email{mazjb@bristol.ac.uk}

\usepackage[T1]{fontenc}
\usepackage[latin1]{inputenc}
\usepackage{latexsym}
\usepackage[arrow, matrix, curve]{xy}
\usepackage[dvips]{graphicx}
\usepackage{epsfig}
\usepackage{bm}
\usepackage{amssymb}
\usepackage{amsmath}
\usepackage{verbatim}
\usepackage{amsthm}
\usepackage{enumerate}

\newtheorem{thm}{Theorem}
\newtheorem{lem}{Lemma}

\newtheorem{cor}{Corollary}

\theoremstyle{remark}
\newtheorem*{rem}{Remark}

\theoremstyle{definition}

\def\B#1{\mathbf{#1}}
\def\ba{\bm{\alpha}}
\def\bb{\bm{\beta}}

\def\bpsi{\bm{\psi}}
\def\F#1{\mathfrak{#1}}
\def\C#1{\mathcal{#1}}
\def\D{\mathrm{d}}
\def\dsum#1#2{\sum_{\substack{{#1}\\{#2}}}}
\def\ol#1{\overline{\B{#1}}}
\def\hs{\B{\hat{h}}}
\def\U#1{\underline{\B{#1}}}

\DeclareMathOperator*{\lcm}{lcm}
\DeclareMathOperator{\card}{Card}
\DeclareMathOperator{\rank}{rank}
\DeclareMathOperator{\vol}{vol}

\newenvironment{pf}{\begin{proof}[Proof]}{\end{proof}}

\usepackage{hyperref}

\begin{document}

\maketitle

\begin{abstract}
    For a given set of forms $\psi^{(1)}, \ldots, \psi^{(R)} \in \mathbb Z[t_1, \ldots, t_m]$ of degree $d$ we prove a Hasse principle for representations of the shape
    \begin{equation*}
        F^{(\rho)}(t_1 \B x_1 + \ldots + t_m \B x_m) = \psi^{(\rho)}(t_1, \ldots, t_m), \quad 1 \leq \rho \leq R
    \end{equation*}
    by general forms $F^{(1)}, \ldots, F^{(R)}  \in \mathbb Z[x_1, \ldots, x_s]$ of the same degree, provided that $s \gg R^2m^d$ and the forms $F^{(\rho)}$ are `sufficiently non-singular'.
    This result is then used to derive asymptotical behaviour of the number of $m$-dimensional linear spaces contained in the intersection of the $F^{(\rho)}$ if the degree is odd. A further application dispenses with the non-singularity condition and establishes the existence of $m$-dimensional linear spaces on the intersection of $R$ cubic forms if the number $s$ of variables asymptotically exceeds $R^6 + m^3 R^3$. Finally, we briefly consider linear spaces on small systems of quintic equations.

\end{abstract}

\section{Introduction}
The problem of determining whether two given forms represent one another is a
classical one and has triggered important developments in the history of modern
number theory. In fact, Gauss's theory of binary quadratic forms can be interpreted
in this way, as can Lagrange's Theorem and related results. Another instance, more
important for our purposes, is the wide range of questions connected to Waring's
problem, which have led to the development of new and powerful sets of tools.

In the 1930s Siegel wrote a series of papers \cite{siegel1,siegel2,siegel3} to
deduce what is now known as Siegel's mass formula, which gives a general description
for indefinite quadratic forms. This has been supplemented by
the discussion of the definite case by Ellenberg and Venkatesh \cite{ellenberg},
based on ideas from ergodic theory, so that the quadratic case is now relatively well understood. Less is known in higher-degree situations, apparently mainly due to the fact that
the methods used in the quadratic case cease to be applicable as soon as the degree
exceeds two. Existing results include the work of Parsell \cite{sp-smooth}, who
examines the possibility of representing a given form of degree $d$ as a sum of
$d$-th powers, using ideas from Arkhipov and Karatsuba \cite{arx}.

In a more general setting over $\mathbb C$, this problem can be addressed with methods from algebraic geometry and has almost completely been solved by Alexander and Hirschowitz \cite{alg-geom-1}. In particular, they proved that with a small number of known and well-understood exceptions, every homogeneous polynomial $\psi(t_1, \ldots, t_m)$ of degree $d$ that can be related to a set of points in general position (see \cite{miranda} for details) has the expected number of representations as a sum of $s$ powers of linear polynomials, provided that
\begin{equation*}
    s \geq  \frac{1}{m} \binom{m-1+d}{d} ,
 \end{equation*}
and this bound is sharp.
Less is known for those forms that do not fulfil the stated generality condition, and finding a lower bound for $s$ that applies to both general and exceptional polynomials is still an unsolved problem even in the complex setting (see the discussion in the introduction of \cite{alg-geom-2}, for instance).
For our purposes, these results are of interest inasmuch as they purvey inherited upper bounds for the number of representations that hold over $\mathbb Q$ as a subfield of $\mathbb C$. However, since the particular structure of $\mathbb Q$ is forfeited by the embedding into the complex numbers, the Alexander-Hirschowitz Theorem is unfit to deliver any real number-theoretic information.\\

In this paper we aim to fill this gap by establishing rather general conclusions
for forms that are, in a certain sense, not too singular. In particular, we will
make use of the Hardy-Littlewood circle method to prove an asymptotic for the number
of identical representations
\begin{equation}\label{rep}
    F^{(\rho)}(t_1 \B x_1  + \ldots + t_m \B x_m ) = \psi^{(\rho)}(t_1, \ldots, t_m) \quad (1 \leq \rho \leq R)
\end{equation}
of a given set of forms $\psi^{(\rho)} \in \mathbb Z[t_1, \ldots, t_m]$ of
degree $d$ by forms $F^{(\rho)} \in \mathbb Z[x_1, \ldots, x_s]$ in $s$ variables.
In order to give a rigorous enunciation of the result we need to introduce some notation.
Let $P$ be a large positive integer, write $\bpsi$ for the $R$-tuple
$\left(\psi^{(1)}, \ldots, \psi^{(R)}\right)$ and denote by $N_{s,R,m}^{(d)}(P; \B F; \bpsi)$
the number of integral solutions of the equations \eqref{rep} with $\B x_i \in
[-P,P]^s$ for $1 \leq i \leq m$. For the sake of convenience we will in the future suppress most of the parameters and use the more concise notation $N_{s,\bpsi}(P)$ or alternatively $N_{s,R,m}^{(d)}(P)$ in the case $\bpsi = \bm 0$.
Furthermore, let $V$ denote the intersection of the singular loci $V^{(\rho)}$ of
the forms $F^{(\rho)}$, and write
\begin{equation} \label{r}
        r = \binom{d-1+m}{d}
\end{equation}
for the number of coefficients of each of the $\psi^{(\rho)}$. Notice that $r \sim m^d$ as $m$ tends to infinity. In this notation the above statement can be quantified in the following way.

\begin{thm} \label{Hasse}
    Let $d \geq 3$, $R$, and $m \geq 2$ be positive integers, and let
        \begin{equation*}
        s-\dim V > 3 \cdot 2^{d-1}  (d-1)R(Rr+1).
    \end{equation*}
     Then there exist nonnegative constants $\chi_{\infty}(\bpsi)$ and $\chi_p(\bpsi)$ for every prime $p$ such that
    \begin{equation*}
        N_{s,\bpsi}(P) = P^{ms-Rrd} \chi_{\infty}(\bpsi) \prod_{p \text{ prime}} \chi_p(\bpsi) + o\left(P^{ms-Rrd}\right).
    \end{equation*}
\end{thm}

The proof of Theorem \ref{Hasse} goes along well-trodden paths that have been
paved in the early sixties by the classical works of Davenport
\cite{dav32,dav29,dav16} and Birch \cite{birch} and, some twenty years later,
developed further by Schmidt \cite{schmidt85}. Note that by expanding and equating powers of the $t_i$, each equation in the system \eqref{rep} takes the shape of a system of $r$ equations. The exponent is therefore the expected one, with $ms$ being the total number of variables and $Rrd$ the total degree of the expanded system of equations.
The critical observation in our argument turns out to be the fact that the system \eqref{rep} can be read either as characterising linear $m$-spaces on $R$ equations or as describing point solutions to an expanded system of $Rr$ equations, and whereas most authors hitherto followed the latter interpretation when implementing the circle method, we will switch freely between the two. This allows us to exploit the fact that the expanded system associated to an equation that describes a linear space will be equipped with a particular structure, which saves considerably over the simpler approach that neglects this structural information by treating the system as if all equations were genuinely distinct.  \\

%
An obvious application of Theorem \ref{Hasse} is that of counting linear spaces on
hypersurfaces. This problem has received considerable attention since the seminal
work of Birch \cite{birch57}, in which he applied a diagonalisation method to prove
the existence of arbitrarily many hyperplanes of any given dimension on the
intersection of an arbitrary number of hypersurfaces, provided only that the number
of variables be sufficiently large, and in order to steer clear of obstructions to
the real solubility, he requires the degrees of the hypersurfaces to be odd. This
method is extremely wasteful in the number of variables (a quantified version has
been provided by Wooley \cite{tdw-birch}), but whereas improved results have been obtained by specialising to simpler special cases (e.g. \cite{rd-tdw}), the understanding of the general situation remains unsatisfactory and consequently few attempts have been made to go beyond the mere existence of linear $m$-spaces and find quantitative estimates.

%
\begin{thm}\label{asymp}
    Let $F^{(1)}, \ldots, F^{(R)} \in \mathbb{Z}[x_1, \ldots, x_s]$ be forms of equal odd degree $d$, and $m \geq 2$ an integer, and let $r$ be as in \eqref{r}. Furthermore, suppose that
            \begin{equation}\label{as-bd}
                s - \dim V > \begin{cases}3 \cdot 2^{d-1} (d-1)  d^{2^{d-1}} R \left(R^2d^2+Rm\right)^{2^{d-2}} & \hbox{ for } d \geq 5, \\ 24R \max \left \{ 10(6R^2+mR)^2, Rr+1 \right\} & \hbox{ if }d = 3. \end{cases}
            \end{equation}
     Then we have
            \begin{equation} \label{asymp-bound}
                N_{s,R,m}^{(d)}(P) =  P^{ms-Rrd} \chi_{\infty} \prod_{p \text{ prime}} \chi_p + o\left(P^{ms-Rrd}\right),
            \end{equation}
    and the product of the local densities $\chi_{\infty} \prod_{p} \chi_p$ associated to $\bpsi = \bm 0$ is positive.\\
 \end{thm}
The case distinction in the cubic case arises from the somewhat surprising fact that the special geometry of our problem can be exploited in a way that enables us to ensure the existence of local $p$-adic solutions for some choices of $m$ and $R$ with much looser conditions on the number of variables than what is needed to establish a Hasse principle, where $s$ is needed to grow at least cubically in $m$ by Theorem \ref{Hasse} and \eqref{r}.
Furthermore, note that due to the generality of the setting, the bound given in \eqref{as-bd} will not be sharp for any typical set of parameters and in any given special case the numerical constants can be improved just by inserting the available bounds for the respective situations.

In order to put this result into context, the most relevant seems to be the work of Dietmann \cite{rd10}, who succeeded in showing that the number of variables necessary in order to guarantee the existence of an affine $m$-space on a single form $F$ grows polynomially in $m$, provided that $F$ is non-singular. Apart from imposing a looser nonsingularity condition, our Theorem \ref{asymp} supersedes Dietmann's bound of
\begin{equation*}
    s \geq 2^{5+2^{d-1}d}d!d^{2^d+1}m^{d(1+2^{d-1})}
\end{equation*}
by a power of $2d$ in $m$, due mainly to our more careful perusal of Schmidt's methods \cite{schmidt85}. \\



The situation is significantly more complicated if one tries to obtain unconditional results.
Denote by $\gamma_d(R, m)$ the least integer $\gamma$ such that any set of $R$ forms $ F^{(\rho)} \in \mathbb Z[x_1, \ldots, x_s]$ of equal degree $d$ in $s>\gamma$ variables contains a rational linear space of affine dimension $m$. Considering the cubic case, Lewis and Schulze-Pillot \cite{lsp} proved that
\begin{equation}\label{lsp-bound}
    \gamma_3(R,m) \ll R^{11}m+R^3m^5,
\end{equation}
while a more geometric approach by Schmidt \cite{lsp} establishes
\begin{equation}\label{schmidt-bound}
\gamma_3(R,m) \ll R^5m^{14}.
\end{equation}
In 1997, Wooley \cite{tdw-book} was able to obtain $\gamma_3(R,m) \ll
R^{8+\epsilon}m^5$, or alternatively $\gamma_3(R,m) \ll_R m^{\alpha}$ with $\alpha=(5+\sqrt{17})/2 = 4.56155...$, and most recently Dietmann \cite{rd08} established a Hasse principle for the number of linear spaces on hypersurfaces which enabled him to show that
\begin{equation}\label{rd-bound}
    \gamma_3(R,m) \ll R^6m^5+R^4m^6.
\end{equation}
We are going to refine Dietmann's methods and prove the following.
\begin{thm}\label{cubic}
    One has
    \begin{equation*}
        \gamma_3(R,m) \ll R^6 + R^3 m^3.
    \end{equation*}
\end{thm}
This is stronger than both \eqref{lsp-bound} and \eqref{rd-bound} and supersedes \eqref{schmidt-bound} for all
$m \gg R^{1/14}$. Furthermore, it will be clear from the proof that the constant can with little effort be made explicit. Note that the bound in Theorem \ref{cubic} is of cubic growth in $m$, thus missing the expected true growth rate only by a power of one. While quadratic growth has been established for the problem concerning a single form in \cite{rd-tdw}, Theorem \ref{cubic} is new even in the case $R=2$, superseding a previous result by Wooley who bounded $\gamma_3(2,m)$ by a quartic polynomial (\cite{tdw97}; see also the discussion of the corollary in \cite{tdw-book}). \\

Following the example of Dietmann \cite{rd10}, we can use Theorems \ref{asymp} and \ref{cubic} and apply an iterating argument to derive an unconditional bound for a quintic version of Theorem \ref{cubic}.
\begin{thm}\label{quintic}
  We have
        \begin{equation*}
            \gamma_5(R,m) \ll_R  m^{12(3^{R-1}-1)+48 \cdot 3^{R-1}}.
        \end{equation*}
  In particular, $\gamma_5(1,m) \ll m^{48}.$
\end{thm}
This is the first time that polynomial growth in $m$ has been established for the problem of finding $m$-spaces on systems of quintic forms, thus improving Dietmann's treatment \cite{rd10} of the case $R=1$, for which he requires at least $\gamma_5(1,m) \ll m^{439}$ variables, both in quantity and in quality. In fact, by a more careful analysis it is possible track the dependence on $R$, and the same methods will yield a bound of the general shape
    \[\gamma_5(R,m) \leq (A mR)^{B ^{c R}}\]
with explicit constants $A,B,c$. For large $R$ this bound is not very satisfactory, as the advantage stemming from polynomial behaviour in $m$ will soon be nullified by the number $R$ of equations occurring in the second order exponent. However, in the light of the work by Wooley \cite{tdw-birch}, especially the discussion in section 6 of his paper, this drawback comes as no surprise. In fact, Wooley's bound \cite{tdw-book} of the shape
\begin{equation*}
    \gamma_5(R,m) \ll (3mR)^{A(mR)^{c}}
\end{equation*}
with numerical constants $A$ and $c$ will prevail as soon as $R \gg \log m$.\\

I should like to thank my supervisor Trevor Wooley for suggesting this problem to
me. Without his constant encouragement and many enlightening discussions this work
would not have been possible. I am also very grateful to Rainer Dietmann for useful conversations on the topics presented in this paper and related ideas and in particular for making available to me a preprint of \cite{rd-neu}. Finally, I would like to express my gratitude to the referee for valuable comments.

\section{Notation and Setting} \label{notation}
Throughout this paper the following notational conventions will be observed: Expressions like $\sum_{n=1}^x f(n)$, where $x$ may or may not be an integer, are always understood to mean $\sum_{1 \leq n \leq x} f(n)$. Also, all estimates including an $\epsilon$ are true for any $\epsilon > 0$, so the same symbol will be applied in all instances. Furthermore, inequalities involving vectors are to be interpreted componentwise, i.e. $ \left| \B a \right| \leq P$ means $\left| a_i \right| \leq P$ for all $i$, and $\gcd(\B a, q)$ is $\gcd(a_1, \ldots, a_n, q)$. We are confident that no misunderstandings should arise if similar statements be read in a like manner. Finally, the exponential $e(x)$ denotes $e(x)=e^{2 \pi i x}$, and the Landau and Vinogradov symbols $O, o, \sim, \gg, \ll, \asymp$ will be used with their established meanings.\\

For $1 \leq \rho \leq R$ let $F^{(\rho)} \in \mathbb Z [x_1, \ldots, x_s]$ be given by
\begin{equation*}
    F^{(\rho)}(\B x)=\sum_{\B i \in \{1, \ldots,s\}^d} c^{(\rho)}(\B i) x_{i_1} \cdot \ldots \cdot  x_{i_d}
\end{equation*}
with symmetric coefficients $c^{(\rho)}(\B i) \in \mathbb Z / d!$, and define the multilinear form $\Phi^{(\rho)}$ associated to $F^{(\rho)}$ by
\begin{equation*}
    \Phi^{(\rho)}\big(\B x^{(1)}, \ldots ,\B x^{(d)}\big) = \sum_{\B i \in \{1, \ldots, s\}^d} c^{(\rho)}(\B i) x_{i_1}^{(1)} \cdot \ldots \cdot  x_{i_d}^{(d)}.
\end{equation*}
Thus one has
\begin{equation*}
    F^{(\rho)}(\B x) = \Phi^{(\rho)}(\B x, \ldots, \B x).
\end{equation*}

In order to count solutions to \eqref{rep}, one needs to understand expressions of the shape
\begin{equation}\label{linspace}
    F^{(\rho)}(t_1 \B x_1+ \ldots + t_m \B x_m).
\end{equation}
This requires an appropriate kind of index notation. Write $J$ for the set of multi-indices $(j_1, j_2, \ldots, j_d) \in \{1, 2, \ldots, m\}^d$, where we allow repetitions in the tuples $(j_1, j_2, \ldots, j_d)$ but disregard order; the number of these is $\card(J)=r$, which is the parameter defined in \eqref{r}. By means of the Multinomial Theorem, \eqref{linspace} can be written as
\begin{equation}\label{expanded}
    F^{(\rho)}\left(t_1 \B x_1+ \ldots + t_m \B x_m\right)
    =\sum_{ \B j \in J} A(\B j) t_{j_1}t_{j_2} \cdot \ldots \cdot t_{j_d} \Phi^{(\rho)}(\B x_{j_1}, \B x_{j_2},\ldots, \B x_{j_d}),
\end{equation}
where the factors $A(\B j)$ take account of the multiplicity of each term and are defined as follows. To every $\B j \in J$ one can associate numbers $\mu_1(\B j), \ldots, \mu_m(\B j)$ between $0$ and $d$ such that
\begin{equation}\label{def-mu}
    t_{j_1}t_{j_2} \ldots t_{j_d} = t_1^{\mu_1(\B j)} t_2^{\mu_2(\B j)} \cdot \ldots \cdot t_m^{\mu_m(\B j)}.
\end{equation}
In other words, the $\mu_i(\B j)$ count the multiplicity with which any given $\B x_i$ appears in the term with index $\B j$.
In this notation, the factors $A(\B j)$ are given by the multinomial coefficients
\begin{equation*}
    A(\B j) = \binom{d}{\mu_1(\B j), \mu_2(\B j), \ldots, \mu_m(\B j)}.
\end{equation*}

Let $\psi^{(1)}, \psi^{(2)}, \ldots, \psi^{(R)}$ be homogeneous polynomials of degree $d$ in $m$ variables, defined by
\begin{equation*}
    \psi^{(\rho)}(t_1, \ldots, t_m)=\sum_{\B j \in J} n^{(\rho)}_{\B j} A(\B j)t_{j_1} t_{j_2} \cdot \ldots \cdot t_{j_d}
\end{equation*}
with integral coefficients $n_{\B j}^{(\rho)}$, and write
\begin{equation*}
    \bpsi=\big(\psi^{(1)}, \psi^{(2)}, \ldots, \psi^{(R)}\big).
\end{equation*}
We can now expand \eqref{rep} and sort by coefficients. This yields a system of equations
\begin{equation*}
    \Phi^{(\rho)}(\B x_{j_1}, \ldots, \B x_{j_d}) = n_{\B j}^{(\rho)} \quad (\B j \in J, 1 \leq \rho \leq R),
\end{equation*}
which is amenable to a circle method approach, so for given $\rho$ write $\ba^{(\rho)}=(\alpha_{\B j}^{(\rho)})_{\B j \in J}$ and let
\begin{align}\label{def-F}
    \F F^{(\rho)}\big(\B x_1, \ldots, \B x_m; \ba^{(\rho)}\big)= \sum_{ \B j \in J} \alpha^{(\rho)}_{\B j} \Phi^{(\rho)}(\B x_{j_1}, \ldots, \B x_{j_d}).
\end{align}
Summing over the $\rho$ is still going to yield a lengthy formula, which can be abbreviated to
\begin{equation}\label{F-final}
    \F F\left(\ol x; \U{\ba}\right) = \sum_{\rho=1}^R \F F^{(\rho)}\big(\B x_1, \ldots, \B x_m; \ba^{(\rho)}\big),
\end{equation}
where we introduced the shorthand notation
    \[ \left(\B x_1, \ldots, \B x_m\right) =\ol x \]
and
    \[ \big(\ba^{(1)}, \ba^{(2)}, \ldots, \ba^{(R)}\big) =\U{\ba}, \]
respectively.
Since it will be useful at some point to sort the components of $\U{\ba}$ by the $\B j$ as opposed to the $\rho$, we seize the opportunity to define $\U \alpha_{\B j}= \big(\alpha^{(1)}_{\B j}, \ldots, \alpha^{(R)}_{\B j}\big)$ for all $\B j \in J$. The same notational conventions will be observed for the coefficients $n_{\B j}^{(\rho)}$ of the target polynomials $\bpsi$.

The expression in \eqref{F-final} collects all the $Rr$ terms that arise from expanding each of the $R$ equations as a sum of $r$ multilinear forms, and thus allows us to define the exponential sum in a very compact notation as
\begin{equation*}
    T(\U{\ba}) = \sum_ {\ol x}  e(\F F \left(\ol x; \U{\ba}\right)).
\end{equation*}
In general, the sum will be over a box $-P \leq \B x_i \leq P$ for all $1 \leq i \leq m$, but in special cases we will write $T(\U{\ba}, X)$ or $T(\U{\ba}, \F B)$ to denote the sum over an $ms$-dimensional hypercube with sidelength $2X$ or a domain $\F B \subset \mathbb Z^{ms}$, respectively.
Altogether, classical orthogonality relations imply that the number of simultaneous representations of $\psi^{(\rho)}$ by $F^{(\rho)}$ contained in the hypercube $[-P,P]^{ms}$ is described by the integral
\begin{align}\label{numbersol}
    N_{s,\bpsi}(P) &= \int_{[0,1)^{Rr}} T(\U{\ba}) e(-\U{\ba} \cdot \U n)\D  \U{\ba} \nonumber \\
    &= \dsum{\B x_1, \ldots, \B x_m}{\left| \B x_i \right| \leq P}  \int_{[0,1)^{Rr}}e\left( \F F\left(\ol x; \U{\ba}\right)-\U{\ba} \cdot \U n\right) \D  \U{\ba}.
\end{align}


It should be noted, however, that although expressions as in \eqref{def-F} and \eqref{F-final} aim to simultaneously solving $rR$ equations, the single equations can be reassembled and can thus be read in the way of our original problem of finding $m$-dimensional linear spaces on the intersection of $R$ hypersurfaces. In fact, one can view the coefficients $\alpha_{\B j}^{(\rho)}$ as absorbing the factors $A(\B j)t_{j_1}\cdot\ldots\cdot t_{j_d}$ arising in \eqref{expanded} and write, somewhat imprecisely,
\begin{align*}
    \F F^{(\rho)}\big(\ol x; \ba^{(\rho)}\big) - \ba^{(\rho)} \cdot \B n^{(\rho)}  &= \sum_{ \B j \in J} \alpha^{(\rho)}_{\B j} \big( \Phi^{(\rho)}(\B x_{j_1}, \ldots, \B x_{j_d}) - n_{\B j}^{(\rho)}\big) \\
    & = \alpha^{(\rho)} \left( F^{(\rho)}(t_1 \B x_1+ \ldots + t_m \B x_m) - \psi^{(\rho)}(t_1, \ldots, t_m) \right).
\end{align*}
Thus \eqref{numbersol} can be read either as finding simultaneous solutions of $Rr$ equations or equivalently as counting $m$-dimensional linear spaces on $R$ hypersurfaces, and while for the greater part of the analysis we will stick to the former interpretation, it will be convenient to switch to the latter one when analysing the singular series more carefully. \\

\section{Weyl differencing}\label{sec-Weyl}
The proof of Theorem \ref{Hasse} is largely along the lines of the classical arguments of Birch \cite{birch} and Schmidt \cite{schmidt85}, with most of the analysis and notation following Birch, while imitating Schmidt's arguments in the treatment of the singular series and singular integral.

The first step is to establish an inequality of Weyl type as presented in chapters 12 -- 14 of \cite{dav}, or in a more general version, in \cite{birch}. Although this is fairly standard, we will give a rather detailed exposition, because it is here that the specific shape of the forms assembled in $\F F(\ol x; \U{\ba})$ comes into play.
\begin{lem} \label{diff1}
    Let $1 \leq k \leq d$ and $j_l$ $(l=1, \ldots, k)$ be integers with $1 \leq j_l \leq m$. Then
    \begin{equation*}
        |T(\U{\ba})|^{2^k} \ll P^{\left((2^k-1)m-k\right)s} \sum_{\B h_1, \ldots, \B h_k \in [-P,P]^s}  \sum_{ \ol x  } e(\Delta_{j_k, \B h_k} \ldots \Delta_{j_1, \B h_1} \F F(\ol x; \U{\ba})),
    \end{equation*}
    where the discrete differencing operator $\Delta_{i, \B{h}}$ is defined by its action on the form $\F F(\ol x; \U{\ba})$ as
    \begin{equation}\label{diff-op}
        \Delta_{i, \B h} \F F(\ol x; \U{\ba}) = \F F(\B x_1, \ldots, \B x_{i} + \B{h}, \ldots, \B x_m; \U{\ba}) - \F F(\B x_1, \ldots, \B x_{i}, \ldots, \B x_m;\U{\ba})
    \end{equation}
    and the sum over the $\ol x$ extends over suitable boxes of sidelength at most $2P$.
\end{lem}
\begin{pf}
    As is usual with Weyl differencing arguments, we proceed by induction. The case $k=1$ follows from a simple application of Cauchy's inequality and one has
    \begin{align*}
        |T(\U{\ba})|^2 \ll& P^{(m-1)s} \dsum{\B x_i }{i\neq j_1}  \left( \sum_{\B h_1} \sum_{\B x_{j_1}}  e\left(\Delta_{j_1, \B h_1} \F F(\ol x; \U{\ba})\right) \right)
    \end{align*}
    by the linearity of $\Delta_{i, \B h}$. Here the summation over $\B{x}_{j_1}$ extends over the intersection of the boxes $|\B{x}_{j_1}| \leq P$ and $|\B{x}_{j_1}+ \B h_{1}| \leq P$, which is again a box of sidelength at most $2P$.

    Now let us assume that the lemma is true for a given $k$. Again by Cauchy's inequality, one finds
    \begin{align*}
        |T(\U{\ba})|^{2^{k+1}} &\ll P^{2\left((2^k-1)m-k\right)s} \left| \sum_{\B h_1, \ldots, \B h_k} \sum_{\ol x} e\left(\Delta_{j_k, \B h_k} \ldots \Delta_{j_1, \B h_1} \F F(\ol x; \U{\ba})\right) \right|^2 \\
        &\ll P^{\left(2^{k+1}-2\right)ms-2ks}  P^{ ks + (m-1)s} \\
        & \qquad \times \sum_{\B h_1, \ldots, \B h_k} \dsum{\B x_i }{i \neq j_{k+1}} \left|\sum_{\B x_{j_{k+1}} }e \left(\Delta_{j_k, \B h_k} \ldots \Delta_{j_1, \B h_1} \F F(\ol x; \U{\ba}) \right)\right|^2 \\
        &\ll P^{\left((2^{k+1}-1)m-(k+1)\right)s} \\
        & \qquad \times  \sum_{\B h_1, \ldots, \B h_{k+1}} \sum_{\ol x}e \left(\Delta_{j_{k+1}, \B h_{k+1}} \ldots \Delta_{j_1, \B h_1} \F F(\ol x; \U{\ba}) \right)
    \end{align*}
    as required.
 \end{pf}
For the sake of notational brevity in the following considerations, we will write $\hs$ for the $(d-1)$-tuple $(\B h_1, \ldots, \B h_{d-1})$. Our final estimate of the exponential sum $T(\U{\ba})$ is an application of the above.
\begin{lem}\label{diff2}
    For any $\B j \in J$ one has the estimate
    \begin{equation*}
        |T(\U{\ba})|^{2^{d-1}} \ll P^{\left(2^{d-1}m-d\right)s} \sum_{\hs} \prod_{i=1}^s \min \left( P, \bigg\| M(\B j) \sum_{\rho=1}^R \alpha_{\B j}^{(\rho)} B_i^{(\rho)}(\hs)\bigg\|^{-1} \right) ,
    \end{equation*}
    where the functions $B^{(\rho)}_i$ are given by
    \begin{equation*}
        \Phi^{(\rho)}(\B x, \B h_1, \ldots, \B h_{d-1}) = \sum_{i=1}^s x_i B_i^{(\rho)}(\B h_1, \ldots, \B h_{d-1})
    \end{equation*}
    and the coefficients $M(\B j)$ are defined by means of \eqref{def-mu} as
    \begin{equation*}
         M(\B j) = \mu_1(\B j)! \mu_2(\B j)! \cdot \ldots \cdot \mu_m(\B j)!.
    \end{equation*}
\end{lem}
\begin{pf}
    Inserting $k=d-1$ in the above lemma gives
    \begin{equation*}
        |T(\U{\ba})|^{2^{d-1}} \ll P^{\left((2^{d-1}-1)m-(d-1)\right)s}  \sum_{\hs} \sum_{\ol x} e\left(\Delta_{j_{d-1}, \B h_{d-1}} \ldots \Delta_{j_1, \B h_1} \F F(\ol x; \U{\ba})\right) .
    \end{equation*}
    By the definition \eqref{diff-op} of $\Delta_{j,\B h}$ and the polynomial structure of $\F F$ every differencing step reduces the degree of the resulting form by one, and therefore this last expression depends only linearly on the vectors $\B x_1, \ldots, \B x_m$. In particular, for any given $\rho$ all the forms assembled in $\F F^{(\rho)}(\ol x; \ba^{(\rho)})$ are instances of the same multilinear form $\Phi^{(\rho)}$ associated to the original form $F^{(\rho)}$, and this structure is naturally preserved by the differencing procedure. Writing $R(\hs)$ for the terms independent of $\ol x$, one obtains
    \begin{align*}
        & \sum_{ \hs } \sum_{\ol x} e\left( \Delta_{j_{d-1}, \B h_{d-1}} \ldots \Delta_{j_1, \B h_1} \F F(\ol x; \U{\ba})\right)\\
        =& \sum_{\hs}\dsum{\B x_i}{i \neq j_d} \sum_{\B x_{j_d}} e \left( M(\B j) \sum_{\rho = 1}^R \sum_{k=1}^m \alpha_{(j_1, \ldots, j_{d-1}, k)}^{(\rho)}\Phi^{(\rho)}(\B x_{k}, \B h_1, \ldots, \B h_{d-1}) + R(\hs)  \right) \\
        \ll& P^{s(m-1)} \sum_{\hs} \left| \sum_{\B x_{j_d}} e \left( M(\B j) \sum_{\rho = 1}^R \alpha_{\B j}^{(\rho)}\Phi^{(\rho)}(\B x_{j_d}, \B h_1, \ldots, \B h_{d-1})  \right) \right|.
    \end{align*}
    Since the $\B x_{j}$ are contained in suitable boxes contained in $[-P,P]^s$, by standard arguments one arrives at the estimate
    \begin{align*}
         |T(\U{\ba})|^{2^{d-1}} &\ll P^{\left(2^{d-1}m-d\right)s} \sum_{\hs} \prod_{i=1}^s \min\left(P, \bigg\| M(\B j) \sum_{\rho=1}^R \alpha_{\B j}^{(\rho)} B_i^{(\rho)}(\hs) \bigg\|^{-1} \right),
    \end{align*}
    which is indeed the required expression.
 \end{pf}

Note that whereas the exponential sum $T(\U{\ba})$ is independent of~$\B j$, this is not true for the expression on the right hand side of the equation. This is due to the fact that the Weyl differencing argument forces us to close in onto one index $\B j$, but since it is immaterial which difference is taken in each step, the estimate holds for any index $\B j \in J$. This pays off when we read the estimate in Lemma \ref{diff2} from the right, as it will save us the averaging over the indices $\B j$, and one obtains independent estimates of equal quality for all $\B j \in J$. This is the main point in which the linear space setting behaves differently from the the general situation as treated in \cite{birch}, as it allows us to treat what is technically a system of $r$ distinct equations as a single one, and this different behaviour will ultimately yield the improvements this paper obtains over previous estimates.

We use Lemma \ref{diff2} to generate a tripartite case distinction. Let us assume that we have $|T(\U{\ba})| \gg P^{ms-k\theta}$ for some parameters $k, \theta>0$, and let $\B j$ be fixed. In this case Lemma \ref{diff2} yields
\begin{align*}
     \left(P^{ms-k\theta}\right)^{2^{d-1}} &\ll P^{(2^{d-1}m-d)s} \sum_{\hs} \prod_{i=1}^s \min\left(P, \bigg\| M(\B j) \sum_{\rho=1}^R \alpha_{\B j}^{(\rho)} B_i^{(\rho)}(\hs) \bigg\|^{-1} \right) ,
\end{align*}
which is equivalent to
\begin{equation*}
    \sum_{\B h_1, \ldots, \B h_{d-1}} \prod_{i=1}^s \min\left(P, \bigg\| M(\B j) \sum_{\rho=1}^R \alpha_{\B j}^{(\rho)} B_i^{(\rho)}(\hs) \bigg\|^{-1} \right) \gg P^{ds-2^{d-1}k\theta}.
\end{equation*}
Note that this expression is independent of $m$. This implies that the following arguments can be extracted directly from standard references. In particular, one can apply the geometry of numbers in order to count how often the minimum takes on a nontrivial value (see \cite{dav}, chapters 12 and 13 for a detailed exposition). The subsequent lemma is the analogue of Lemma 13.3 of \cite{dav} or Lemma 2.4 of \cite{birch}, respectively.

\begin{lem}\label{geom}
    Suppose that $|T(\U{\ba})| \gg P^{ms-k\theta}$ for some parameters $k, \theta>0$, and let $N(X,Y)$ denote the number of $(d-1)$-tuples $\B h_1, \ldots, \B h_{d-1}$ in the box $|\B h_k| \leq X$ satisfying
    \begin{equation}\label{dioph. approx.}
        \bigg\| M(\B j) \sum_{\rho=1}^R \alpha_{\B j}^{(\rho)} B_i^{(\rho)}(\hs) \bigg\| < Y
    \end{equation}
    for all $i=1, \ldots, s$ and $\B j \in J$. Then we have the estimate
    \begin{equation*}
        N(P^{\theta}, P^{-d+(d-1)\theta}) \gg P^{(d-1)s\theta-2^{d-1}k\theta - \epsilon}.
    \end{equation*}
\end{lem}

From Lemma \ref{geom} one infers that either the exponential sum is small, or for all $i$ and $\B j$ the quantity $ M(\B j) \sum_{\rho=1}^R \alpha_{\B j}^{(\rho)} B_i^{(\rho)}(\hs)$ is often close to an integer. The latter can be effected in two ways, as it will occur either if the forms $B_i^{(\rho)}$ tend to vanish for geometric reasons, or by genuine (i. e. non-zero) solutions to the diophantine approximation problem that is implicit in \eqref{dioph. approx.}. This yields a threefold case distinction, which is the heartpiece of all circle method arguments concerning general homogeneous polynomials.
\begin{lem}\label{ram1}
    Let $0 < \theta \leq 1$ and $k>0$ be parameters, and let $\U{\ba} \in [0,1)^{rR}$. Then there are three possibilities:
    \begin{enumerate}[(A)]
    \item
        The exponential sum $T(\U{\ba})$ is bounded by
        \begin{equation*}
            |T(\U{\ba})| \ll P^{ms-k\theta}.
        \end{equation*}

    \item
        For every $\B j \in J$ one finds $(q_{\B j}, \underline{a}_{\B j}) \in \mathbb Z^{R+1}$ satisfying
        \begin{align*}
            q_{\B j} &\ll P^{(d-1)R\theta}  \quad \hbox{and}\\
            \left|\alpha_{\B j}^{(\rho)} q_{\B j} - a_{\B j}^{(\rho)}\right| &\ll P^{-d+(d-1)R\theta} \quad \text{ for all }1 \leq \rho \leq R.
        \end{align*}

    \item \label{sing}
        The number of $(d-1)$-tuples $(\B h_1, \ldots, \B h_{d-1}) \leq P^{\theta}$ that satisfy
        \begin{equation} \label{rank}
            \rank \left(B_i^{(\rho)}(\B h_1, \ldots, \B h_{d-1})\right)_{i, \rho} \leq R-1
        \end{equation}
        is asymptotically greater than $\left(P^{\theta}\right)^{(d-1)s - 2^{d-1}k - \epsilon}$.
    \end{enumerate}
\end{lem}

\begin{pf}
    This follows by the same argument as Lemma 2.5 in \cite{birch}. Suppose that the estimate in (A) does not hold, so that by Lemma \ref{geom} for every $\B j \in J$ we have
    \begin{equation*}
        \bigg| M(\B j) \sum_{\rho=1}^R \alpha_{\B j}^{(\rho)} B_i^{(\rho)}(\hs_{\B j}) \bigg| < P^{-d+(d-1)\theta}
    \end{equation*}
    for at least $\gg P^{(d-1)s\theta-2^{d-1}k\theta - \epsilon}$ choices of $\hs \leq P^{\theta}$. Fixing some index $\B j \in J$ and writing $B(\B j)$ for the $(R \times s)$-matrix $\left(M(\B j)B_i^{(\rho)}(\hs_{\B j})\right)_{i, \rho}$, this means we can find integer vectors $\B A(\B j)$ and $\bm{\delta}(\B j) \in \mathbb Z^s$ with the property that
    \begin{equation}\label{approx1}
         B(\B j) \U{\alpha}_{\B j}   - \B A({\B j})  = \bm{\delta}(\B j)\ll P^{-d+(d-1)\theta}.
    \end{equation}

    If the matrix $B(\B j)$ is of full rank for some $(d-1)$-tuple $\hs_{\B j} = (\B h_1, \ldots, \B h_{d-1})$, then we can find a non-vanishing $(R \times R)$-minor whose absolute value we denote by $q_{\B j}$. We remark here for further reference the obvious fact that $q_{\B j}$ is independent of any particular index $\rho$. This allows us to implement a discrete matrix inversion in order to generate approximations of the~$\U{\alpha}_{\B j}$.

    We may assume without loss of generality that the non-vanishing minor $B_0(\B j)$ of $B(\B j)$ is the first one, and write $\B A_0(\B j), \bm{\delta}_0(\B j) \in \mathbb Z^R$ for the corresponding portions of $\B A(\B j)$ and $\bm{\delta}(\B j)$. Then we can find integer solutions $a_{\B j}^{(1)}, \ldots, a_{\B j}^{(R)}$ to the system
    \begin{equation}\label{approx2}
          B_0(\B j)\underline a_{\B j} = q_{\B j}\B  A_0(\B j).
    \end{equation}

    Combining \eqref{approx1} and \eqref{approx2}, we obtain
    \begin{equation*}
           B_0(\B j)\left(q_{\B j}\U{\alpha}_{\B j}   - \underline a_{\B j} \right) =  q_{\B j}  \bm{\delta}_0(\B j).
    \end{equation*}
    By Cramer's rule this returns the required bound, and the proof is complete on noting that
    \begin{equation*}
        q_{\B j} \ll \max_{i, \rho} |B_i^{(\rho)}(\hs_{\B j})|^R \ll P^{R(d-1)\theta}.
    \end{equation*}

    Note that since the estimates obtained in Lemma \ref{diff2} hold independently of the index $\B j$, the differencing variables $\hs_{\B j}$ that generate the rational approximation need not be the same for all $\B j \in J$, and therefore the approximations will in general depend on the $\B j$-component of the $\alpha_{\B j}^{(\rho)}$ that is being approximated. On the other hand, it is clear that the bounds themselves, in particular the upper bound for the values of $B_i^{(\rho)}(\hs_{\B j})$, are independent of the index $\B j$ chosen in the beginning of the proof, so we can find rational approximations of the same quality for all vectors $\U{\alpha}_{\B j}$, and their denominators vary with $\B j \in J$ but are independent of $\rho$.
 \end{pf}

The condition \eqref{rank} of case \eqref{sing} in Lemma \ref{ram1} is tantamount to a system of simultaneous equations in $s(d-1)$ variables and thus defines a variety which we call $\F V$. By \cite{birch}, Lemma 3.2, its dimension is
\begin{equation}\label{schmidt-g}
    \dim \F V \geq (d-1)s  - 2^{d-1}k .
\end{equation}
Since, however, the variety $\F V$  is not particularly easy to handle, we replace it by the intersection of the singular loci $V$ as in \cite{birch}, Lemma 3.3, and obtain
\begin{equation*}
    s-\dim V \leq (d-1)s - \dim \F V \leq 2^{d-1}k .
\end{equation*}
This allows us to exclude the third case in Lemma \ref{ram1} by choosing the number of variables sufficiently large.
\begin{lem} \label{ram2}
    Let $\U{\ba} \in [0,1)^{Rr}$ and let $0 < \theta \leq 1$ and $k>0$ be parameters with
        \begin{equation}\label{k-cond}
            s-\dim V > 2^{d-1}k.
        \end{equation}
    Then the alternatives are the following:
    \begin{enumerate} [(A)]

    \item \label{minor}
        The exponential sum $T(\U{\ba})$ is bounded by
        $$|T(\U{\ba})| \ll P^{ms-k\theta}.$$

    \item  \label{major}
        For every $\B j \in J$ one finds $(q_{\B j}, \underline{a}_{\B j}) \in \mathbb Z^{R+1}$, satisfying
        \begin{align*}
            q_{\B j} &\ll P^{(d-1)R\theta}  \quad \hbox{and} \\
            \left|\alpha_{\B j}^{(\rho)} q_{\B j} - a_{\B j}^{(\rho)}\right| &\ll P^{-d+(d-1)R\theta} \quad (1 \leq \rho \leq R).
        \end{align*}

    \end{enumerate}
\end{lem}
%
%

\section{Major Arcs dissection} \label{major arcs dissection}

%

Lemma \ref{ram2} suggests a major arcs dissection in terms of the parameter $\theta$, a notion that can be made rigorous by specifying the implicit constant. Let $C$ be sufficiently large in terms of the coefficients of the $F^{(\rho)}$. We define the major arcs $\F M(P,\theta)$ to be the set of all $\U{\ba} \in [0,1)^{Rr}$ that have a rational approximation satisfying
\begin{align}\label{approx}
    0 \leq a_{\B j}^{(\rho)} <q_{\B j} &\leq C P^{(d-1)R\theta}  \quad \hbox{and} \nonumber \\
    \left|\alpha_{\B j}^{(\rho)} q_{\B j} - a_{\B j}^{(\rho)}\right| &\leq CP^{-d+(d-1)R\theta} \quad (1 \leq \rho \leq R),
\end{align}
and the minor arcs
\begin{equation*}
    \F m(P,\theta) = [0,1)^{Rr} \setminus \F M(P,\theta)
\end{equation*}
to be the complement thereof. In order to save clout, we will omit the parameter $P$ during most of the analysis, specifying it only in cases where ambiguities might be likely to arise. It is, however, worthwhile to note that this definition respects the case distinction of Lemma \ref{ram2}, that is, for every $\U{\ba} \in [0,1)^{rR}$ one has either a rational approximation as in \eqref{approx} or the the estimate in case \eqref{minor} holds true.
By the argument of \cite{birch}, Lemma 4.1, the major arcs are disjoint if $2R(d-1)\theta <d$ and their volume is at most

\begin{align} \label{volmajor}
    \mathrm{vol} (\F M(\theta)) &\ll  \prod_{{\B j} \in J} \left(\sum_{q_{\B j}=1}^{CP^{R(d-1)\theta}} \prod_{\rho=1}^R \left(\sum_{a_{\B j}^{(\rho)}=0}^{q_{\B j}-1} \frac{P^{-d+R(d-1)\theta}}{q_{\B j}} \right) \right) \nonumber\\
    &\ll \prod_{{\B j} \in J} \left(\sum_{q_{\B j}=1}^{CP^{R(d-1)\theta}} \left(P^{-d+R(d-1)\theta}\right)^R \right) \nonumber\\
    &\ll P^{-drR+(d-1)rR(1+R)\theta}.
\end{align}

As will become apparent in the following discussion, we will need to fix the parameter $\theta$ rather small so as to allow a better error control when examining the major arcs contribution more closely. Also, in order to minimise the number of variables required in \eqref{k-cond}, we should like to choose $k$ small. On the other hand, we require $k \theta > Rrd$ in order to get a suitable estimate on the minor arcs. This discrepancy motivates the following pruning lemma.

\begin{lem} \label{pruning}
    Suppose the parameters $k$ and $\theta$ satisfy
    \begin{align*}
        0 < \theta & < \theta_0 = \frac{d}{(d-1)(R+1)}
    \end{align*}
    and
    \begin{equation}
        k >  Rr(R+1)(d-1). \label{k}
    \end{equation}
    Then there exists a $\delta > 0$ such that the minor arcs contribution is bounded by
    \begin{equation*}
        \int_{\F m(P,\theta)} |T(\U{\ba})| \D \U{\ba} \ll P^{ms-Rrd-\delta}.
    \end{equation*}
\end{lem}
\begin{pf}
    This is a straightforward adaptation of Lemma 4.4 in \cite{birch}.
    Given $0 < \theta < \theta_0$, we can find a parameter $\delta >0$ such that
    \begin{equation} \label{delta1}
        (k-Rr(R+1)(d-1))\theta > 2\delta
    \end{equation}
    and a sequence $\theta_i$ with the property that
    \begin{equation*}
        1 \geq \theta_0 > \theta_1 > \theta_2 > \ldots > \theta_{M-1} > \theta_M = \theta >0
    \end{equation*}
    and subject to the condition
    \begin{equation} \label{delta2}
        (\theta_i-\theta_{i+1})k<\delta \quad \text{for all $i$}.
    \end{equation}
    This is always possible with
    \begin{equation} \label{M}
        M=O(1).
    \end{equation}
    Then on writing
    \begin{equation*}
        \F m_i = \F m(\theta_i) \setminus \F m(\theta_{i-1}) = \F M(\theta_{i-1}) \setminus \F M(\theta_i)
    \end{equation*}
    one has
    \begin{equation*}
        \mathrm{vol} (\F m_i)  \leq \mathrm{vol} (\F M(\theta_{i-1})) \ll P^{-drR+(d-1)rR(1+R)\theta_{i-1}}
    \end{equation*}
    by \eqref{volmajor}. Recall that for $\U{\ba} \in \F m(\theta)$, we are in the situation of case (A) in Lemma \ref{ram2}, so the minor arcs contribution is bounded by
   \begin{align*}
        \int_{\F m(\theta) \setminus \F m(\theta_0)} \left| T(\U{\ba}) \right| \D \U{\ba} &= \sum_{i=1}^M \int_{\F m_i} \left| T(\U{\ba}) \right| \D \U{\ba}  \\
        &\ll \sum_{i=1}^M \mathrm{vol} (\F M(\theta_{i-1})) \sup_{\U{\ba} \in \F m(\theta_i) } \left| T(\U{\ba}) \right| \\
        &\ll \sum_{i=1}^M P^{-Rrd+(d-1)rR(1+R)\theta_{i-1}} P^{ms-k\theta_i}.
   \end{align*}
   By \eqref{M}, the sum is of no consequence and can be replaced by a maximum over all $i \in \{1, \ldots, M\}$. Hence the exponent is
   \begin{align*}
        & - Rrd+(d-1)rR(1+R) \theta_{i-1} + ms-k\theta_i\\
        = & \, ms-Rrd+k(\theta_{i-1}-\theta_i) -(k-(d-1)rR(1+R))\theta_{i-1}\\
        \leq & \, ms-Rrd-\delta,
   \end{align*}
    where the last inequality uses \eqref{delta1} and \eqref{delta2}. The result follows on noting that
    \begin{align*}
        \int_{\F m(\theta_0)} \left| T(\U{\ba}) \right| \D \U{\ba} \ll P^{ms-k \theta_0} = o(P^{ms-Rrd}).
   \end{align*}
\end{pf}

\section{Homogenising the approximations}\label{homogenising}

Lemma \ref{ram2} \eqref{major} gives us approximations of the shape
\begin{equation*}
    \alpha_{\B j}^{(\rho)} = a^{(\rho)}_{\B j}/q_{\B j} + \beta_{\B j}^{(\rho)}
\end{equation*}
with denominators that are in general different for each ${\B j} \in J$. It will, however, greatly facilitate the future analysis if we can find a common denominator $q$ such that approximations of the shape
\begin{equation*}
    \alpha_{\B j}^{(\rho)} = b_{\B j}^{(\rho)}/q + \gamma_{\B j}^{(\rho)}
\end{equation*}
and of a similar quality hold. For sufficiently small $\theta$ this is indeed possible, but in order to homogenise the set of major arcs we have to surmount some technical difficulties.

Define
\begin{align}\label{b-def}
    q = \lcm_{\B j \in J} \{ q_{\B j} \} \quad \hbox{ and } \quad  b_{\B j}^{(\rho)} = a_{\B j}^{(\rho)} q / q_{\B j} \quad (\B j \in J, 1 \leq \rho \leq R),
\end{align}
and note that $\gcd(\U b, q)=1$.

\begin{lem}\label{weight}
	Let $q$ and $b_{\B j}^{(\rho)}$ be as above. There exist integer weights $\lambda_{\B j}^{(\rho)}$ for all $\B j \in J$ and $1 \leq \rho \leq R$ such that $\lambda_{\B j}^{(\rho)} \leq q_{\B j}$ and
    \begin{equation}\label{gcd}
        \gcd\Big(\sum_{\B j, \rho}  \lambda_{\B j}^{(\rho)} b_{\B j}^{(\rho)}, q\Big)= 1.
    \end{equation}
\end{lem}

\begin{pf}
    By Euclid's algorithm there exist parameters $\lambda_{\B j}^{(\rho)}$ such that
        $$\sum_{\B j, \rho} \lambda_{\B j}^{(\rho)} b_{\B j}^{(\rho)} = \gcd_{\B j, \rho}\left\lbrace  b_{\B j}^{(\rho)} \right\rbrace $$
    and hence
        $$\gcd\Big( \sum_{\B j, \rho}  \lambda_{\B j}^{(\rho)} b_{\B j}^{(\rho)},q\Big) = \gcd(\U b,q) =1. $$
    These $\lambda_{\B j}^{(\rho)}$ really live modulo $q_{\B j}$, since by writing $\lambda_{\B j}^{(\rho)} = c_{\B j}^{(\rho)} q_{\B j} + \mu_{\B j}^{(\rho)}$ and recalling \eqref{b-def} one has
        $$ b_{\B j}^{(\rho)} \lambda_{\B j}^{(\rho)} = \frac{a_{\B j}^{(\rho)} q}{ q_{\B j}} \left(c_{\B j}^{(\rho)} q_{\B j} + \mu_{\B j}^{(\rho)}\right) \equiv \frac{a_{\B j}^{(\rho)} q}{ q_{\B j}}\mu_{\B j}^{(\rho)} \equiv b_{\B j}^{(\rho)}\mu_{\B j}^{(\rho)} \pmod q.$$
    This allows us to take $\lambda_{\B j}^{(\rho)} \leq q_{\B j}$.
 \end{pf}

The following Lemma may be useful in other contexts, so we will state it in a rather general fashion.

\begin{lem}\label{CH'}
    Let $\F B $ be the image of $[-P,P]^n \cap \mathbb Z^n$ under some integral non-singular linear transformation $A$, and assume that $\F B  \subset [-X,X]^{n}$ for some $X$. Furthermore, let
        $$T(\ba,\F B) = \sum_{\B x \in \F B} e(\ba \cdot \B F(\B x) )$$
    be a multidimensional exponential sum over $\F B $.
    Then we have
    \begin{equation*}
        \big| T(\ba, \F B) \big| \ll \det(A)^{-1} (X/P)^n ( \log P)^n \sup_{\bm{\eta} \in [0,1)^{n}} \big| H(\ba, \bm{\eta}; X) \big|,
    \end{equation*}
    where $H(\ba, \bm{\eta}; X)$ is defined as
    \begin{equation*}
        H(\ba, \bm{\eta}; X) = \sum_{|\B x| \leq X} e\big( \ba \cdot \B F(\B x) - \bm{\eta} \cdot \B x \big).
    \end{equation*}
\end{lem}
The important fact to notice here is that $H(\ba,\bm{\eta}; P)$ is just the usual exponential sum $T(\ba, P)$ with a linear twist characterised by the parameter $\bm{\eta} \in [0,1)^{n}$. Thus its behaviour will not essentially differ from that of the usual exponential sum, and the two can be regarded as roughly the same object. This means that Lemma \ref{CH'} enables us to treat exponential sums over rather more general convex sets than standard rectangular boxes.

\begin{pf}
    By the orthogonality relations, one has
    \begin{align*}
       \sum_{\B y \in \F B}  e\left(\ba \cdot \B F(\B y) \right)
       &= \sum_{|\B x| \leq X}e\left(\ba \cdot \B F(\B x)\right) \sum_{\B y \in \F B} \int_{[0,1)^{n}} e(\bm{\eta} \cdot \B y - \bm{\eta} \cdot \B x) \D \bm{\eta}\\
       &= \int_{[0,1)^{n}}\sum_{|\B x| \leq X}e\left(\ba \cdot \B F(\B x) -\bm{\eta} \cdot \B x \right) \sum_{\B y \in \F B}  e(  {\bm{\eta}} \cdot \B y) \D \bm{\eta}.
    \end{align*}
    Thus if we write
    \begin{equation*}
        D(\bm{\eta}, \F B)=\sum_{\B y \in \F B}  e( \bm{\eta} \cdot \B y),
    \end{equation*}
    the exponential sum $T(\ba, \F B)$ can be expressed in terms of $H(\ba, \bm{\eta}; X)$ as
    \begin{align*}
        T(\ba, \F B) & = \int_{[0,1)^{n}} H(\ba, \bm{\eta}; X) D(\bm{\eta}, \F B) \D \bm{\eta} \\
        & \ll \sup_{\bm{\eta} \in [0,1)^{n}} \big| H(\ba, \bm{\eta}; X) \big| \int_{[0,1)^{n}} D(\bm{\eta}, \F B) \D \bm{\eta}.
    \end{align*}
    Now the $\B y$ are in the image of $[-P,P]^{n}$ under $A$ and can therefore be written as $\B y = A \B x$ with $| \B x | \leq P$.
    This implies that
        $$\bm{\eta} \cdot \B y = \bm{\eta} \cdot A\B x = A^t\bm{\eta} \cdot \B x$$
    and consequently
    \begin{align*}
        D(\bm{\eta}, \F B) =  D(A^t\bm{\eta}, [-P,P]^n) \ll  \prod_{i=1}^n  \min (P, \| ( A^t \bm{\eta})_i \|^{-1} ).
    \end{align*}
    It follows that
    \begin{align*}
         \int_{[0,1)^{n}} D(\bm{\eta}, \F B) \D \bm{\eta} &\ll \int_{[0,1)^{n}} \prod_{i=1}^n  \min (P, \| ( A^t \bm{\eta})_i \|^{-1} ) \D \bm{\eta} \\
         & \ll (\det A)^{-1} \int_{\F C} \prod_{i=1}^n  \min (P, \|  \eta_i \|^{-1})  \D \bm{\eta},
    \end{align*}
    where $\F C$ is the image of $[0,1)^{n}$ under $A^t$. Obviously, the integrand is positive and $1$-periodic in every direction, so we can bound the integral over $\F C$ by a number copies of the integral over the unit cube, where the factor is determined by the number of unit cubes needed to cover $\F C$. Since $A$ has the property of mapping $[-P,P]^n$ into a subset of $[-X,X]^n$, this is bounded by $(X/P)^n$, and one has
    \begin{align*}
        \int_{\F C} \prod_{i=1}^n  \min (P, \|   \eta_i \|^{-1} ) \D \bm{\eta}  &\ll (X/P)^n \int_{[0,1)^n} \prod_{i=1}^n  \min (P, \|   \eta_i \|^{-1} ) \D \bm{\eta} \\
        &\ll  \left(X/P\right)^n(\log P)^n.
    \end{align*}
    This gives the result.
 \end{pf}

\begin{rem}
    Note that in the last step one has $\vol(\F C) = \det(A)$, so by periodicity one would really expect something like
    \begin{equation*}
        \int_{\F C} \prod_{i=1}^n  \min (P, \|   \eta_i \|^{-1})  \D \bm{\eta} \approx \det(A) \int_{[0,1)^n} \prod_{i=1}^n  \min (P, \|   \eta_i \|^{-1} ) \D \bm{\eta}
    \end{equation*}
    to hold, which would yield
    \begin{equation*}
        \big| T(\ba, \F B) \big| \ll (\log P)^n \sup_{\bm{\eta} \in [0,1)^{n}} \big| H(\ba, \bm{\eta}; X) \big|
    \end{equation*}
    in the statement of the lemma. Unfortunately, this heuristic relies on the assumption that the integrand is well-behaved over general sets, but since we are not making any further hypotheses regarding the transformation $A$, it cannot be taken for granted that the fluctuations of the integrand cancel out. For instance, it is possible for the image of the unit cube under $A$ to be stretched along one of the axes, so that it contains more than the expected number of near-integer points and thus gives a much greater contribution than expected. It is in order to account for this that one needs the worse multiplicity factor $(X/P)^n$ instead of $\det (A)$.
\end{rem}

We have now collected the technical tools necessary in order to homogenise our set of major arcs, which allows us to proceed and prove a homogenised version of Lemma \ref{ram2}.

\begin{lem}\label{hom-lemma}
   Assume that $\theta R(d-1)(r+3)<d.$ Under the conditions of Lemma \ref{ram2}, we can replace alternative \eqref{major} by
   \begin{enumerate}
   \item[(B$'$)]
        There exists an integer $q \ll P^{2(d-1)R\theta}$ such that one finds $\U{\B a} \in \mathbb Z^{Rr}$, satisfying
        \begin{align*}
            \left|\alpha_{\B j}^{(\rho)} q - a_{\B j}^{(\rho)}\right| &\ll P^{-d+3(d-1)R\theta} \quad (1 \leq \rho \leq R, \, \B j \in J).
        \end{align*}
   \end{enumerate}
\end{lem}

%


%

\begin{pf}
    Make the substitution
    \begin{align}\label{def-x'}
        \B x_k' = \begin{cases} \B x_k-\B x_{k+1} & 1 \leq k \leq m-1 \\ \B x_m & k=m, \end{cases}
    \end{align}
    so that
    \begin{align*}
        \B x_k &= \sum_{i=k}^m \B x_i' \quad (1 \leq k \leq m).
    \end{align*}
    Furthermore, observe that the index set $J$ is equipped with a partial order relation. We may order the entries of the multi-indices by size, i.e. $j_1 \leq j_2 \leq \ldots \leq j_d$, then the partial order is induced by entrywise comparison, so we say that $\B j \leq \B j'$ for two elements $\B j, \B j' \in J$ if and only if $j_k \leq j_k'$ for all $1 \leq k \leq d$.
    
    Since the proof of Lemma \ref{ram1} produces the same denominators $q_{\B j}$ independently of the $\rho$-component of the 	coefficients $\alpha_{\B j}^{(\rho)}$, it suffices without loss of generality to consider only the case $\rho=1$. This allows us to avoid unnecessary complexity of the notation by dropping the index and writing $\ba$ instead of $\ba^{(1)}$.
    
    Observe that
    \begin{align*}
        \Phi(\B x_{j_1}, \ldots, \B x_{j_d}) = \Phi\Big( \sum_{k_1 \geq j_1} \B x'_{k_1}, \ldots,  \sum_{k_d \geq j_d}\B x'_{k_d}\Big)  =  \sum_{\B k \geq  \B j}  \Phi(\B x'_{k_1}, \ldots, \B x'_{k_d}).
    \end{align*}
	
    Now consider the weighted exponential sum
	\begin{equation*}
        T( \bm{\lambda}\ba, P )= \dsum{|\B x_i| \leq P}{1 \leq i \leq m} e\left( \sum_{\B j \in J} \lambda_{\B j} \alpha_{\B j}\Phi (\B x_{j_1}, \ldots, \B x_{j_d}) \right),
    \end{equation*}
	where the $\lambda_{\B j}$ will be fixed later. This can be expressed in terms of the alternative variables $ \ol x'$ and yields
    \begin{align*}
        T( \bm{\lambda}\ba, P ) = \dsum{\B x_i'}{1 \leq i \leq m} e\left( \sum_{\B k \in J} \Big( \sum_{\B j \leq \B k}  \lambda_{\B j} \alpha_{\B j} \Big) \Phi (\B x'_{k_1}, \ldots, \B x'_{k_d})\right)
    \end{align*}
    where the sum over $\B x_i'$ is over domains $\C B_i$ contained in $[-2P,2P]^{s}$ as determined by \eqref{def-x'}. Thus if we define $\ba'$ by
    \begin{equation*}
        \alpha_{\B j}' = \sum_{\B k \leq \B j} \lambda_{\B k} \alpha_{\B k} ,
    \end{equation*}
    this yields the identity
        $$T( \bm{\lambda}\ba, P ) = T({\ba}',  \textstyle{\prod_i}\C  B_i).$$
    Beware that the domains $\mathcal B_i$ are not independent of one another, so one has to apply great care in exchanging the order of summation, and this affects our possibilities of applying Weyl's inequality severely.
    However, since the transformation \eqref{def-x'} is non-singular and maps $[-P,P]^{ms}$ onto a subset of $[-2P,2P]^{ms}$, Lemma \ref{CH'} comes to our rescue and yields
     \begin{equation}\label{CH-eq}
        |T( \bm{\lambda}\ba, P )| = |T({\ba}',  \textstyle{\prod_i}\C  B_i)| \ll  \displaystyle{(\log P)^{ms} \sup_{\overline{\bm{\eta}} \in [0,1)^{ms}} \big| H({ \ba}', \ol{\bm{\eta}}; 2P) \big|}.
    \end{equation}

    The exponential sum $T({\ba}',P)$ and its twisted cousin $H({\ba}', \ol{\bm{\eta}}, 2P)$ should be thought of as being roughly of the same order of magnitude. Indeed, since
    \begin{equation*}
            \Delta_{j_2, \B h_2}\Delta_{j_1, \B h_1} (\ol{\bm{\eta}} \cdot \ol x) = \Delta_{j_2, \B h_2} (\bm{\eta}_{j_1} \cdot \B h_1) = 0,
    \end{equation*}
    any linear twist has no effect in the deduction of Weyl's inequality, and a small modification of the proof of Lemma \ref{diff1} readily shows that Lemma \ref{diff2} continues to hold if $T({\ba}',P)$ is replaced by $H({\ba}', \ol{\bm{\eta}}, P)$ as long as $d \geq 2$. This implies that all following estimates of section \ref{sec-Weyl} will remain unaffected by the twist. In particular, the minor arcs estimate will hold for $H({\ba}', \ol{\bm{\eta}}, P)$ if and only if it does so for $T({\ba}',P)$. Collecting these arguments together, one concludes that $|T({\ba}',P)| \ll P^{ms-k\theta}$ implies that $|H({\ba}', \ol{\bm{\eta}}, P)| \ll P^{ms-k\theta}$ for arbitrary $\ol{\bm{\eta}}$, and by \eqref{CH-eq} it follows that $|T( \bm{\lambda}\ba, P )| \ll P^{ms-k\theta}(\log P)^{ms}$. Since the singular case is excluded, this in turn means that whenever $ \bm{\lambda}\ba$ possesses an approximation as in case \eqref{major} of Lemma \ref{ram2}, then so does $\ba'$. \\

	Now suppose that $\ba \in \F M(\theta)$ for some $\theta$. Then it has approximations $\alpha_{\B j} = a_{\B j}/q_{\B j} + \beta_{\B j}$, and according to Lemma \ref{weight} we can find integer weights $\lambda_{\B j} \leq q_{\B j}\ll P^{(d-1)R\theta}$ satisfying \eqref{gcd}. One has
    \begin{equation*}
        \lambda_{\B j}\alpha_{\B j}=\lambda_{\B j} \frac{a_{\B j}}{q_{\B j}} + \lambda_{\B j}\beta_{\B j},
    \end{equation*}
    so $\bm{\lambda}\ba$ is certainly contained in $\F M(2\theta)$. By the above considerations this implies that we have rational approximations for the components $\alpha_{\B j}'$ of $\ba'$,  given by $a'_{\B j}$ and $q'_{\B j} \ll P^{2(d-1)R\theta}$ such that $|\alpha'_{\B j}q'_{\B j} - a'_{\B j}| \ll P^{-d+2(d-1)R\theta}$.

    Consider the last term $\alpha_{\B m}' = \alpha_{m, \ldots, m}'$.  We have the approximation
    \begin{equation}\label{a' reg-approx}
         \alpha_{\B m}'=  \frac{a'_{\B m}}{q_{\B m}'}+ O\left(  \frac{P^{-d+2(d-1)R\theta}}{q_{\B m}'} \right).
    \end{equation}
	On the other hand, $\B m \geq \B j $ for all $\B j \in J$ with respect to the partial order of the $\B j$, so inserting the definition of $\alpha_{\B m}'$ gives the alternative approximation
	\begin{align}\label{a' alt-approx}
        \alpha_{\B m}' &= \sum_{\B j \in J}  \alpha_{\B j} \lambda_{\B j} = \sum_{\B j \in J} \left( \frac{\lambda_{\B j}a_{\B j}}{q_{\B j}} + O \left(\frac{P^{-d+2R(d-1)\theta}}{q_{\B j}}\right) \right) \nonumber \\
        & = \frac{\sum_{\B j \in J}  \lambda_{\B j} b_{\B j}}{q}+ O\left( \frac{P^{-d+(r+1)(d-1)R\theta}}{q}\right).
    \end{align}
    Our goal is now to show that the two approximations \eqref{a' reg-approx} and \eqref{a' alt-approx} are actually the same. However, if they are distinct, it follows that
    \begin{align*}
         \frac{1}{q_{\B m}' q} & \leq \left|\frac{a_{\B m}'}{q_{\B m}'} - \frac{\sum_{\B j \in J}  \lambda_{\B j} b_{\B j}}{q} \right|  \\
         &\leq \left| \alpha_{\B m}'- \frac{a_{\B m}'}{q_{\B m}'} \right| + \left| \alpha_{\B m}'  -\frac{\sum_{\B j \in J}  	\lambda_{\B j} b_{\B j}}{q} \right| \\
         &\ll \frac{qP^{-d+2(d-1)R\theta} + q_{\B m}'P^{-d+(r+1)(d-1)R\theta}} {q_{\B m}'q}\\
         & \ll \frac{P^{-d+(r+3)(d-1)R\theta}} {q_{\B m}'q} ,
    \end{align*}
    which is a contradiction if $\theta$ is sufficiently small. Choosing $\theta$ in accordance with the hypothesis of the statement of the lemma, we can thus conclude that \eqref{a' reg-approx} and \eqref{a' alt-approx} coincide. Lemma \ref{weight} now ensures that $\gcd(\sum_{\B j \in J}  \lambda_{\B j} b_{\B j}, q) = 1$, so the above approximations are both reduced fractions and one has $q =  q_{\B m}' \ll P^{2(d-1)R\theta}$.

    Finally, the bound on $|q\alpha_{\B j}- b_{\B j}|$ follows by observing that
    \begin{equation*}
        |\alpha_{\B j}q-b_{\B j}| = \frac{q}{q_{\B j}} |\alpha_{\B j} q_{\B j}-a_{\B j}| \leq q |\alpha_{\B j} q_{\B j}-a_{\B j}| \ll 	P^{2(d-1)R\theta}P^{ -d+ (d-1)R\theta}.
    \end{equation*}
    This yields the statement.
 \end{pf}

It should be noted here that the factors $2$ and $3$ that arise in the homogenising process probably have little right to exist at all, and that one would expect them to succumb to a more momentous argument than the ones we have been presenting here. However, the important accomplishment of this section is to avoid collecting another factor $r$ in the homogenising process, as that would undo the gains from Lemma \ref{ram2} and throw us back into the situation of treating \eqref{rep} as a system in $rR$ variables without regard for symmetries.

\section{Generating functions analysis}\label{generating functions analysis}
Our goal in this section is to show that the major arcs contribution can be interpreted as a product of local densities.
In order to do so, it is somewhat inconvenient that the estimate of $\U{\bb}= |\U{\ba} - \U a/q|$ depends on $q$. We therefore extend the major arcs slightly and define our final choice of major arcs $\F M'(P,\theta)$ to be set of all $\U{\ba} = \U a/q + \U{\bb}$ contained in $[0,1)^{Rr}$ that satisfy
\begin{align} \label{major'}
    | \U{\bb} | &\leq C' P^{-d+3(d-1)R\theta} \nonumber\\
      \U a < q &\leq C' P^{2(d-1)R\theta}
\end{align}
for some suitably large constant $C'$. As before, we may suppress the parameter $P$. Notice also that by Lemma \ref{hom-lemma} this definition comprises the original major arcs as defined in \eqref{approx}, provided $\theta$ is small enough. Henceforth all parameters $\U{\ba}, \U a, q, \U{\bb}$ will be implicitly understood to satisfy the major arcs inequalities as given in \eqref{major'}.

Letting
\begin{align*}
    S_q(\U a) = \sum_{\ol x=1}^{q} e \left( \frac {\F F(\ol x; \U a)}{q}  \right)
\end{align*}
and
\begin{align*}
    v_P( \U{\bb}) = \int_{|\ol{\bm{\xi}}| \leq P} e \left( \F F(\ol{\bm{\xi}};  \U{\bb}) \right) \D  \ol{\bm{\xi}},
\end{align*}
we can replace the exponential sum by an expression that reflects the rational approximation to $\U{\ba}$ and will be easier to handle.

\begin{lem} \label{gen}
    Assume that $\U{\ba} \in \F M'(P, \theta)$. Then there exists an integer vector $(\U a, q)$ such that
    \begin{align*}
         T(\U{\ba})  - q^{-ms} S_q(\U a) v_P( \U{\bb}) &\ll q^{ms} \left( 1 + (Pq^{-1})^{ms} qP^{d-1} | \U{\bb} | \right)  \\
         &\ll P^{ms-1+5(d-1)R\theta} .
    \end{align*}
\end{lem}
\begin{pf}
    The first estimate is essentially like Lemma 8.1 in \cite{sp-smooth} by sorting the variables into arithmetic progressions and applying the Mean Value Theorem, whereas the second inequality follows from inserting the major arcs estimates \eqref{major'} for $\U{\bb}$ and $q$ and noting that $\theta$ is small enough for the first term $q^{ms}$ to be negligible.
 \end{pf}
The next step is to integrate the expression from Lemma \ref{gen} over $\F M'(P,\theta)$ in order to determine the overall error arising from this substitution. For this purpose, define the truncated singular series and singular integral as
\begin{align*}
    \F S_{\bpsi}(P) &= \sum_{q=1}^{C'P^{2(d-1)R\theta}} q^{-ms} \dsum{\U a=0}{(\U a, q)=1}^{q-1} S_q(\U a)  e\left(-(\U{\B n} \cdot \U{\B a})/q \right)
\end{align*}
and
\begin{align*}
    \F J_{\bpsi}(P) &= \int_{|\U{\bb}| \leq C'P^{-d+3(d-1)R\theta}}v_P(\U{\bb}) e(-\U{\B n} \cdot \U{\bb}) \D  \U{\bb},
\end{align*}
respectively.
\begin{lem}\label{generror}
   The total major arcs contribution is given by
    \begin{equation*}
        \int_{\F M'(P,\theta)} T (\U{\ba}) e(-\U{\ba} \cdot \U n) \D  \U{\ba} = \F S_{\bpsi}(P) \F J_{\bpsi}(P) + O\left(P^{ms-Rrd+(d-1)R(5Rr+7)\theta -1}\right).
    \end{equation*}
\end{lem}
The error is acceptable if $\theta$ has been chosen small enough.

\begin{pf}
    The volume of the extended major arcs $\F M'(\theta)$ is bounded by
    \begin{align*}
            \mathrm{vol} (\F M'(\theta)) &\ll  \sum_{q=1}^{C'P^{2R(d-1)\theta}} \prod_{{\B j} \in J} \prod_{\rho=1}^R \left(\sum_{a_{\B j}^{(\rho)}=0}^{q-1} P^{-d+3R(d-1)\theta} \right)  \\
            &\ll  \sum_{q=1}^{C'P^{2R(d-1)\theta}} \left(q P^{-d+3R(d-1)\theta}\right)^{rR}  \\
            &\ll P^{-Rrd+(d-1)R(5Rr+2)\theta}.
    \end{align*}
    This, together with Lemma \ref{gen}, implies the statement.
 \end{pf}

We can now fix $\theta$ such that $(d-1)R(5Rr+7)\theta < 1$ for the rest of our considerations. Note that with this choice Lemmata \ref{pruning}, \ref{hom-lemma} and \ref{generror} are applicable and one has
\begin{align}\label{asymp-formula}
    N_{s,\bpsi}(P) = \F S_{\bpsi}(P) \F J_{\bpsi}(P) + o(P^{ms-Rrd}),
\end{align}
provided the conditions \eqref{k-cond} and \eqref{k} are satisfied.\\

The truncated singular series and integral can be extended to infinity.
Recalling the definition in (\ref{def-F}), a standard computation reveals that
\begin{equation*}
    \F J_{\bpsi}(P) =  P^{ms-Rrd} \int_{|\U{\bb}| \leq C'P^{3(d-1)R\theta}} v_1(\U{\bb}) e(-\U{\B n} \cdot \U{\bb}) \D  \U{\bb}
\end{equation*}
(see \cite{dav}, Lemma 4.3, for instance). We therefore define, if existent, the complete singular series $\F S_{\bpsi}$ and the singular integral $\F J_{\bpsi}$ as
\begin{align*}
     \F S_{\bpsi} &=\sum_{q=1}^{\infty} \dsum{\U a=0}{(\U a, q)=1}^{q-1} q^{-ms}S_q(\U a) e\left(-(\U{\B n} \cdot \U{\B a})/q \right)  \\
     \F J_{\bpsi} &= \int_{\mathbb R^{rR}} v_1(\U{\bb}) e(-\U{\B n} \cdot \U{\bb}) \D  \U{\bb}.
\end{align*}
In either case, convergence implies that the errors $|\F S_{\bpsi}(P)-\F S_{\bpsi}|$ and \linebreak $\left|\F J_{\bpsi}(P)-P^{ms-Rrd}\F J_{\bpsi}\right|$ are $o(1)$, and we will be able to replace the statement in Lemma \ref{generror} by
\begin{equation*}
    \int_{\F M'(P,\theta)} T (\U{\ba}) e(-\U{\ba} \cdot \U n)  \D  \U{\ba} = P^{ms-Rrd}\F J_{\bpsi} \F S_{\bpsi} +o\left(P^{ms-Rrd}\right).
\end{equation*}
It remains to show that the above definitions are permissible.

As a first step in that direction, we note that by standard arguments one has
\begin{equation*}
    S_{q_1}(\U a_1)S_{q_2}(\U a_2) = S_{q_1q_2}(q_2 \U a_1+ q_1 \U a_2)
\end{equation*}
for coprime $q_1$ and $q_2$, so we have the multiplicativity property
\begin{equation}\label{mult}
    \dsum{\U a_1=1}{(\U a_1,q_1)=1}^{q_1} S_{q_1} (\U a_1) \dsum{\U a_2=1}{(\U a_2,q_2)=1}^{q_2} S_{q_2}(\U a_2) = \dsum{\U b=1}{(\U b,q_1q_2)=1}^{q_1q_2}  S_{q_1q_2}(\U b).
\end{equation}
This allows us to restrict ourselves to considering prime powers in the analysis of the singular series.

\begin{lem} \label{SS-lemma}
    Let $d \geq 3$, $p$ prime and $h$ a non-negative integer, and suppose further that \eqref{k-cond} holds true. Then for any $W>0$ such that $k > (d-1)RW$ the terms of the singular series are bounded by
    \begin{equation*}
        (p^h)^{-ms}S_{p^h}(\U a) \ll (p^h)^{-W}.
    \end{equation*}
\end{lem}

\begin{pf}
    We imitate the argument of Lemma 7.1 in \cite{schmidt85}. Pick a suitable $\theta_1 < (R(d-1))^{-1}$ such that $k\theta_1 \geq W$, and assume that the argument $\U a/p^h$ is on the corresponding major arcs $\F M(p^h,\theta_1)$. Then by the definition \eqref{approx} of the major arcs one can find $ q_{\B j} \ll p^{h(d-1)R\theta_1} \ll p^h$ and $\underline b_{\B j} <  q_{\B j}$ for each $\B j \in J$ subject to
    \begin{equation*}
        \left| q_{\B j} \frac{a_{\B j}^{(\rho)}}{p^h}-b_{\B j}^{(\rho)}\right| \ll (p^h)^{-d+(d-1)R\theta_1} \quad ( 1 \leq \rho \leq R).
    \end{equation*}
    Inserting the bound on $\theta_1$ yields
    \begin{equation*}
        |q_{\B j} a_{\B j}^{(\rho)}-b_{\B j}^{(\rho)}p^h| \ll (p^h)^{2-d}
    \end{equation*}
    for all $\B j \in J$ and $1 \leq \rho \leq R$. For $d>2$ this forces $ q_{\B j} a_{\B j}^{(\rho)}$ to be a multiple of $p^h$, but since $q_{\B j} < p^h$ it follows that $p$ divides $a_{\B j}^{(\rho)}$ for every set of indices. This is, however, impossible as we had $\gcd( \U a,p)=1$. By Lemma \ref{ram2} we can therefore conclude that the minor arcs estimate is true and
    \begin{equation*}
        (p^h)^{-ms}S_{p^h}(\U a) \leq (p^h)^{-k\theta_1} \ll (p^h)^{-W}
    \end{equation*}
    as claimed.
 \end{pf}

\begin{lem}\label{conv-sing-int}
    Suppose $k > 3(d-1)RW$ for some $W>r$, and assume \eqref{k-cond} is satisfied. Then
    \begin{equation*}
        v_1(\U{\bb}) \ll \left(1+\|\U{\bb}\| \right)^{-W}.
    \end{equation*}
\end{lem}
\begin{pf}
    Just as in the proof of Lemma 11 in \cite{schmidtquad}, we note that the equation
    \begin{equation*}
        v_1(\U{\bb}) = Q^{-ms}v_Q(Q^{-d}\U{\bb})
    \end{equation*}
    holds for arbitrary $Q$.
    For
    \begin{equation}\label{theta2-cond}
        0<\theta_2\leq \frac{d}{(d-1)R(r+3)}
    \end{equation}
    (to be determined later) and for a given $\U{\bb} \in \mathbb R^{rR}$ choose $Q$ such that for a suitable constant $C$ one has
    \begin{equation} \label{Q=beta}
         C Q^{3(d-1)R\theta_2}=| \U{\bb} |.
    \end{equation}
    Notice that the condition \eqref{theta2-cond} on $\theta_2$ allows us to apply Lemma \ref{hom-lemma}, so we can assume the major arcs to be homogenised. Hence with this choice of $Q$, the argument $Q^{-d}\U{\bb}$ lies just on the edge of the corresponding major arcs $\F M'(Q,\theta_2)$; in fact, it is best approximated by $\U a= \U 0$ and $q=1$, and one has $q^{-ms}S_1(\U{\bm 0})=1$. An application of Lemma \ref{gen} yields
    \begin{equation*}
        v_Q(Q^{-d}\U{\bb}) = T(Q^{-d} \U{\bb}) +O\left( Q^{ms-1+5(d-1)R\theta_2}\right).
    \end{equation*}

    On the other hand, since $Q^{-d}\U{\bb}$ lies on the edge of the extended the major arcs $\F M'(Q, \theta_2)$, it is not contained in the original set $\F M(Q,\theta_2)$ of major arcs. We can therefore bound the exponential sum $ T(Q^{-d} \U{\bb})$ by the minor arcs estimate and find
    \begin{equation*}
        T(Q^{-d} \U{\bb}) \ll Q^{ms-k\theta_2} \ll Q^{ms-3(d-1)RW\theta_2}.
    \end{equation*}
    This gives
    $$v_1(\U{\bb}) \ll Q^{-ms} \left( Q^{ms-3(d-1)RW\theta_2} +Q^{ms-1+5(d-1)R\theta_2}\right),  $$
    which is optimised by picking
    \begin{equation*}
        \theta_2^{-1}=5(d-1)R+3(d-1)RW = R(d-1)(3W+5).
    \end{equation*}
    Notice that this choice satisfies \eqref{theta2-cond} as we assumed $W>r$. This allows us to rewrite \eqref{Q=beta} in the shape
    \begin{equation*}
        C Q^{3/(3W+5)}=| \U{\bb}|,
    \end{equation*}
    whence the bound on $v_1$ is
    \begin{align*}
        v_1(\U{\bb}) &\ll \left( Q^{-3(d-1)RW\theta_2} +Q^{-1+5(d-1)R\theta_2}\right)  \\
        &\ll Q^{-3W/(3W+5)} \ll | \U{\bb} |^{-W}.
    \end{align*}
    Furthermore, one has trivially $v_1(\U{\bb}) \ll 1$, so on taking the maximum one retrieves the statement.
 \end{pf}

\begin{lem}\label{conv}
    Suppose $k > 3(d-1)R(Rr+1)$. Then the singular series $\F S_{\bpsi}$ and the singular integral $\F J_{\bpsi}$ are absolutely convergent, and one has
    \[ \F S_{\bpsi}(P)\F J_{\bpsi}(P) = P^{ms-rRd}  \big( \F S_{\bpsi}\F J_{\bpsi} +o(1) \big).\]
\end{lem}

\begin{pf}
    This is now immediate from Lemmata \ref{SS-lemma} and \ref{conv-sing-int} on choosing $W>Rr+1$ and noting that $3(Rr+1)>r(R+1)$, so with this choice of $k$ Lemma \ref{pruning} will be applicable.
\end{pf}

Now by the multiplicativity property \eqref{mult} the singular series $\F S_{\bpsi}$ can be expanded as an Euler product
\begin{equation*}
    \F S_{\bpsi}=\prod_p \chi_p(\bpsi),
\end{equation*}
where the $p$-adic densities $\chi_p(\bpsi)$ are given by
\begin{align}
    \chi_p(\bpsi) &= \sum_{i=0}^{\infty} p^{-ims} \dsum{\U a=1}{(\U a, p)=1}^{p^i} S_{p^i}(\U a)e\left(-(\U{\B n} \cdot \U{\B a})/p^i \right) \label{p-adic}  \\
     &= \sum_{i=0}^{\infty}p^{-ims}\dsum{\U a =1}{(\U a, p)=1}^{p^i} \sum_{\ol x =1}^{p^i} e ( \F F(\ol x; p^{-i}\U a)-(\U{\B n} \cdot \U{\B a})/p^i). \nonumber
\end{align}
By the discussion in section 3 of \cite{schmidt85}, this can be interpreted as a $p$-adic integral
\begin{equation*}
    \chi_p(\bpsi)= \int_{\mathbb Q_p^{rR}} \int_{\mathbb Z_p^{ms}} e \left( \F F (\ol{\bm{\xi}}; \U{\bm{\eta}}) -\U{\B n} \cdot \U{\bm{\eta}} \right) \D  \ol{\bm{\xi}} \D  \U{\bm{\eta}},
\end{equation*}
which is the exact analogue of \eqref{numbersol} in the $p$-adic numbers $\mathbb Q_p$. Similarly, the singular integral
\[
    \F J_{\bpsi} = \int_{\mathbb R^{rR}} \int_{|\ol{\bm{\xi}}| \leq 1} e \left( \F F(\ol{\bm{\xi}};  \U{\bb}) -\U{\B n} \cdot \U{\bb} \right) \D  \ol{\bm{\xi}}  \D  \U{\bb}
\]
measures solutions in the real unit box and may thus be interpreted as the density of real solutions. Theorem \ref{Hasse} now follows from \eqref{asymp-formula}, Lemma \ref{conv} and \eqref{k-cond}.

\section{Proof of Theorem \ref{asymp}}\label{local}
In order to prove Theorem \ref{asymp}, it remains to analyse under what conditions the singular series $\F S$ corresponding to $\bpsi = \bm 0$ is positive.
Recalling \eqref{p-adic}, a standard argument (see \cite[Lemma 5.2 and Cor.]{dav} for instance) shows that $\chi_p=1+O(p^{-1-\delta})$ and thus
\begin{equation*}
    1/2 \leq \prod_{p>p_0}\chi_p \leq 3/2
\end{equation*}
for for some suitable $p_0$. It is thus sufficient to show that every individual $\chi_p$ is positive.

For fixed $p$ we split the factors $\chi_p$ and write
\begin{equation*}
    \chi_p=\sum_{i=0}^l \dsum{|\U a| < p^i}{(\U a, p)=1} p^{-ims}S_{p^i}(\U a) + \sum_{i=l+1}^{\infty} \dsum{|\U a| < p^i}{(\U a, p)=1} p^{-ims}S_{p^i}(\U a) = I_l + I_{\infty}.
\end{equation*}
Denote by  $\gamma_d^*(R, m)$ the least integer $\gamma$ such that any set of $R$ forms of equal degree $d$ in $s>\gamma$ variables contains a $p$-adic linear space of affine dimension $m$ for all primes $p$.
Choosing $W>\max\{\gamma_d^*(R, m), Rr\}$ in Lemma \ref{SS-lemma} yields
\begin{align}\label{I_infty}
    I_{\infty} &\ll \sum_{i=l+1}^{\infty} \dsum{|\U a| < p^i}{(\U a, p)=1}  p^{-ims}S_{p^i}(\U a) \nonumber \\
    &\ll \sum_{i=l+1}^{\infty} p^{i(Rr-W )} \nonumber\\
    &\ll p^{l(Rr-\gamma_d^*(R, m) - \delta)}
\end{align}
for some $\delta>0$.

On the other hand, by standard transformations (e.g. \cite{schmidt85} eq. (3.2)) one finds
\begin{equation*}
    I_l = p^{l(Rr-ms)}\Gamma(\F F, p^l),
\end{equation*}
where $\Gamma(\F F, p^l)$ denotes the number of solutions $\ol x \pmod{p^l} $ of the simultaneous congruences
\begin{equation*}
    \Phi^{(\rho)}(\B x_{j_1}, \B x_{j_2}, \ldots, \B x_{j_d}) \equiv 0 \pmod {p^l}, \quad 1 \leq \rho \leq R, \, \B j \in J.
\end{equation*}
%

At this stage, it is useful to recall that the exponential sum in question can be understood not only in terms of counting points on a total of $rR$ equations, but has the alternative and, in fact, more accurate interpretation of describing the number of linear $m$-spaces on $R$ equations. It is clear from the discussion in section~\ref{notation} that by taking suitable linear combinations and applying \eqref{expanded}, the definition of $\Gamma(\F F, p^l)$ can be expressed in terms of the original setting as the number of solutions $\ol x \pmod{p^l} $ of
\begin{equation*}
    F^{(\rho)}( \B x_1 t_1 + \ldots + \B x_m t_m) \equiv 0 \pmod {p^l}, \quad 1 \leq \rho \leq R,
\end{equation*}
where the equivalence is understood to hold identically in $t_1, \ldots, t_m$. This allows us to return to the original formulation of the problem in terms of counting linear spaces.
\begin{lem}\label{schmidt}
    We have
    \begin{equation*}
        \Gamma(\F F, p^l) \gg p^{l\left(ms-\gamma_d^*(R, m)\right)}.
    \end{equation*}
\end{lem}
\begin{pf}
    This is analogous to Lemma 2 in \cite{cubIV}. The proof is by a counting argument that remains intact if statements about $p$-adic points are replaced with the respective statements about $p$-adic linear spaces.
 \end{pf}
As a consequence, one finds that
\begin{equation*}
     I_l \gg p^{l\left(Rr-\gamma_d^*(R, m)\right)}
\end{equation*}
for every prime $p$, which in combination with \eqref{I_infty} yields $ \chi_p = I_l+I_{\infty} \gg 1 $ for a suitable $l$, provided that
\begin{equation} \label{k>v}
    k>3(d-1)RW >3(d-1)R\max \{\gamma_d^*(R, m), Rr \}.
\end{equation}

On the other hand, Lemma \ref{conv} implies that we need $k>3(d-1)R(Rr+1)$ in order to ensure that the singular series $\F S = \prod \chi_p \ll 1$ whence by including the estimate \eqref{k-cond}, the number of variables required is given by
\begin{align*}
    s-\dim(V) > 2^{d-1}k > 3 \cdot 2^{d-1}(d-1)R \max\{\gamma_d^*(R,m), Rr+1\}.
\end{align*}
Finally, $\gamma_d^*(R, m)$ can be controlled by inserting bounds from the literature.
\begin{lem}\label{gamma}
    Let $R, m$ and $d\geq 2$ be positive integers. Then
    \begin{equation*}
        \gamma_d^*(R,m) \leq (R^2d^2+mR)^{2^{d-2}}d^{2^{d-1}} .
    \end{equation*}
    If $d=3$, we have the sharper bound
    \begin{equation*}
        \gamma_3^*(R,m) \leq 10(6R^2+mR)^{2}.
    \end{equation*}
\end{lem}
\begin{pf}
    This is on inserting Theorem 1 of \cite{dav-lewis} or the explicit values given in the subsequent remark, respectively, into Theorem 2.4 of \cite{tdw-local}.
 \end{pf}

Regarding the singular series it can be shown from the Borsuk--Ulam Theorem that for forms of odd degree this will always be positive (see \cite{cubIV}, section 2 for details).
This concludes the proof of Theorem \ref{asymp}. \\

Often it is useful to express results of a shape similar to that of Theorem \ref{asymp} in a way that avoids explicit mention of the singularities. Let therefore $h(F)$ denote the least integer $h$ that allows a form $F$ to be written identically as a decomposition
\begin{equation} \label{decomp}
    F(\B x) = G_1(\B x) H_1(\B x) + \ldots + G_h(\B x) H_h(\B x)
\end{equation}
of forms $G_i, H_i$ of degree strictly smaller than $\deg(F)$. For a system of forms one takes the minimum over the forms in the rational pencil of the $F^{(\rho)}$ and defines
\begin{equation*}
    h(\B F) = \min_{\B c } h\left(c^{(1)}F^{(1)} + \ldots + c^{(R)}F^{(R)}\right),
\end{equation*}
where $\B c$ runs over the non-zero elements of $\mathbb Q^R$. This is a geometric invariant commonly called the $h$-invariant, and it allows us to restate Theorem \ref{asymp} as follows.
\begin{cor}\label{cor}
    Let $F^{(1)}, \ldots, F^{(R)}, R$, $d \geq 5$ and $m$ be as in Theorem \ref{asymp}. Then the
number of linear spaces satisfies \eqref{asymp-bound}, provided that
    \begin{equation*}
        h(\B F) \geq (\log 2)^{-d} d!(d-1)\left(3 \cdot 2^{d-1} d^{2^{d-1}} R (R^2d^2+Rm)^{2^{d-2}} +R(R-1)\right).
    \end{equation*}
\end{cor}

This is readily deduced from Theorem \ref{asymp}. In fact, it follows from Propositions III and III$_{\mathbb C}$ as well as the subsequent corollary in \cite{schmidt85} that $h(\B F)$ is bounded above by
    \[ h(\B F) \leq (\log 2)^{-d}d! \big((s(d-1) - \dim (\F V)) + (d-1)R(R+1)\big).  \]
If the system $\B F$ is singular, then one has
    \[ s(d-1) -\dim (\F V) \leq 2^{d-1}k \]
by \eqref{schmidt-g}. Hence for
\begin{align*}
     h(\B F) &> (\log 2)^{-d}d!\big(2^{d-1}k + (d-1)R(R+1)\big) \\
     &> (\log 2)^{-d}d!\big(3 \cdot 2^{d-1}(d-1)R \max \{ \gamma_d^*(R,m), Rr+1 \}+ (d-1)R(R+1)\big),
\end{align*}
the singular case is excluded. Inserting Lemma \ref{gamma} and noting that for $d \geq 5$ the bound stemming from $\gamma_d^*(R,m)$ dominates in the maximum completes the proof.\\

In the quadratic and cubic case one can do better than in the corollary. In fact, Schmidt's work \cite{schmidtquad,cubIV} on the subject has recently received an improvement by Dietmann \cite{rd-neu}, which translates into our case and which we will apply in our derivation of Theorem \ref{cubic}.
Write
\begin{equation*}
    K(R,m) =24R \max \left \{ 10(Rm+6R^2)^2, Rr+1 \right\}
\end{equation*}
for the bound in the cubic case of Theorem \ref{asymp}.
\begin{thm}\label{rd-thm}
    Let $F^{(1)}, \ldots, F^{(R)}$ be cubic forms in $s$ variables such that no form in the rational pencil vanishes on a linear space of codimension less than $ K(R,m)$. Then we have
    \begin{equation*}
        N_{s,R,m}^{(3)}(P) =  P^{ms-3Rr} \chi_{\infty} \prod_{p \text{ prime}} \chi_p + o\left(P^{ms-3Rr}\right),
    \end{equation*}
    and the product of the local densities $\chi_{\infty} \prod_{p} \chi_p$ is positive.
\end{thm}

This is essentially Theorem 2 of \cite{rd-neu}, but the setting is different enough to warrant some further justification. The proof rests on understanding the singular case of Weyl's inequality for points on the intersection on $R$ cubic forms. For a given $R$-tuple $(w_1, \ldots, w_R)$ let $V(w_1, \ldots, w_R; P)$ denote the number of solutions $\B x^{(1)}, \B x^{(2)} \in [-P,P]^s$ of
\begin{equation*}
    \sum_{\rho=1}^R w_{\rho} B_i^{(\rho)}\big(\B x^{(1)}, \B x^{(2)}\big) = 0 \quad (1 \leq i \leq s),
\end{equation*}
where the $B_i^{(\rho)}$ are as in Lemma \ref{diff2}. Then we have the following alternative version of the singular case of Lemma \ref{ram1}.
\begin{lem}\label{rd-weyl}
    Let $d=3$ and suppose we are in the situation of Lemma \ref{ram1}\eqref{sing}. Then there exist integers $w_1, \ldots, w_R$, not all of which are zero, such that
    \begin{equation*}
        V(w_1, \ldots, w_R; P^\theta) \gg (P^{\theta})^{2s-4k - \epsilon}.
    \end{equation*}
\end{lem}
\begin{pf}
    This is Lemma 2 of \cite{rd-neu} in the case $d=3$. Notice that the singularity condition in Lemma \ref{ram1} does not depend on the value of $m$, so the argument in the singular case is identical.
 \end{pf}
\begin{lem}\label{rd-weyl2}
    Suppose that $d=3$ and each form of the rational pencil of the $F^{(1)}, \ldots, F^{(R)}$ has $h$-invariant greater than $4k$. Then the singular case of Lemma \ref{ram1} is excluded.
\end{lem}
\begin{pf}
    This is essentially Lemma 6 of \cite{rd-neu}. Lemma \ref{rd-weyl} states that in the singular case \eqref{sing} of Lemma \ref{ram1} one can find integers $w_1, \ldots, w_R$ such that
    \begin{equation*}
        V(w_1, \ldots, w_R; P^\theta) \gg (P^{\theta})^{2s-4k- \epsilon}.
    \end{equation*}
    For a fixed set of such $w_{\rho}$ consider the cubic form $C = \sum_{\rho} w_{\rho} F^{(\rho)}$ and notice that this implies that $h(C) \geq h(\B F)$. The proof of Lemma 6 of \cite{rd-neu} shows that the bilinear forms $B_i$ associated to $C$ can be expressed in terms of the original bilinear forms $B_i^{(\rho)}$ associated to the forms $F^{(\rho)}$ as
    \begin{equation*}
        B_i\big(\B x^{(1)}, \B x^{(2)}\big)=\sum_{\rho=1}^R w_{\rho} B_i^{(\rho)}\big(\B x^{(1)}, \B x^{(2)}\big)  \quad (1 \leq i \leq s).
    \end{equation*}
    Lemma 5 of \cite{rd-neu} now implies that
    \begin{equation*}
       \card \left\{\B x^{(1)}, \B x^{(2)} \in [-P,P]^s: B_i\big(\B x^{(1)}, \B x^{(2)}\big) = 0 \quad (1 \leq i \leq s) \right\} \ll P^{2s-h(C)},
    \end{equation*}
    whence we have the sequence of inequalities
    \begin{align*}
         (P^{\theta})^{2s-4k - \epsilon} &\ll V(w_1, \ldots, w_R; P^{\theta}) \\
         & \ll \card \left\{\B x^{(1)}, \B x^{(2)} \in [-P^{\theta},P^{\theta}]^s: B_i\big(\B x^{(1)}, \B x^{(2)}\big) = 0 \quad (1 \leq i \leq s) \right\}\\
         &\ll (P^{\theta})^{2s-h(C)}.
    \end{align*}
    This is a contradiction if $h(C) > 4k$.
 \end{pf}
Theorem \ref{rd-thm} follows now from Lemma \ref{rd-weyl2} exactly as Theorem \ref{asymp} does from Lemma \ref{ram2}. In order to guarantee that the singular series is positive, we need to take $k$ according to \eqref{k>v}, and this together with the bound from Lemma \ref{conv} reproduces the bound $h(\B F)>K(R,m)$. Finally, notice that the condition $h(\B F) \leq K(R,m)$ implies that one of the forms in the rational pencil vanishes on a space of codimension at most $h(\B F) \leq K(R,m)$. This completes the proof of Theorem \ref{rd-thm}.

\section{Linear spaces on the intersection of cubic and quintic hypersurfaces}\label{liner spaces}

Let us now turn our attention to the proof of Theorem \ref{cubic}, and suppose that no form in the rational pencil of the $F^{(1)}, \ldots, F^{(R)}$ vanishes on a rational subspace with codimension at most $ K(R,m)$. Then Theorem \ref{rd-thm} gives $\gamma_3(R,m) \leq K(R,m)$, which proves the theorem for this case. Let us therefore suppose that one of the forms in the linear pencil does vanish on a linear space $Y$ with
\begin{equation*}
    \dim Y \geq s - K(R,m);
\end{equation*}
we can assume without loss of generality that this is the form $F^{(R)}$. This allows us to reduce the problem to finding an $m$-dimensional linear space on the intersection of $Y$ with the hyperspaces associated to $F^{(1)}, \ldots, F^{(R-1)}$ at the expense of having to increase the number of variables by $K(R,m)$. Notice that this is asymptotically bounded above and below by
\begin{equation*}
    K(R,m) \asymp R^5+R^3m^2+R^2m^3 \asymp R^5+R^2m^3.
\end{equation*}
This yields the recursion formula
    \begin{equation*}
        \gamma_3(\rho,m) \leq \max\{K(\rho,m), \gamma_3(\rho-1,m) +K(\rho,m)\} = \gamma_3(\rho-1,m) +K(\rho,m).
    \end{equation*}
Iterating the argument gradually reduces the number of forms, and after at most $R$ steps we retrieve the bound
    \begin{equation*}
        \gamma_3(R,m) \leq \sum_{\rho=1}^R K(\rho, m) \ll  R^6 + R^3m^3.
    \end{equation*}
Furthermore, it is clear that the implied constant is absolute and computable; in fact, a rough calculation confirms that one can take
    \begin{equation*}
        \gamma_3(R,m) \leq 2100 \left( (R+1)^6+(R+1)^3(m+1)^3 \right).
    \end{equation*}\\

Theorem \ref{quintic} is proved similarly. Again, we proceed by induction.
For the case $R=1$ we imitate Dietmann \cite{rd10}, using Theorem \ref{asymp} instead of the weaker bounds applied by him. Applying Corollary \ref{cor} with $d=5$, one sees that any single quintic hypersurface $F$ contains a linear $m$-space as soon as
\begin{equation}\label{hbound}
    h(F) \gg (d^2+m)^{2^{d-2}} \asymp m^8 ,
\end{equation}
so we may suppose that $h(F) \leq C_1 m^{8}$ for some constant $C_1$. By the definition of $h(F)$, one can find forms $G_i, H_i$ $(i=1, \ldots, h)$ of degree less than five such that the form $F$ can be written in the shape given in \eqref{decomp}, and the argument of \cite{rd10} allows us to assume without loss of generality that all forms $H_i$ are cubic. Theorem~\ref{cubic} now implies that the intersection of the $H_i$ contains a linear $m$-space if $s \gg h^3m^3+h^6$. Together with \eqref{hbound} this gives the result.

Now consider $\gamma_5(R,m)$ for $R>1$. As before, in the case $h(\B F) \gg_R m^8$ the claim follows from Theorem~\ref{asymp}. Let us therefore suppose that $F^{(R)}$ possesses a decomposition as in \eqref{decomp} with $h(\B F) \leq C_2(R) m^8$ for some $C_2(R)$.
As above, it suffices to consider the worst case scenario that all forms $H_i^{(R)}$ are cubic, and by Theorem~\ref{cubic} we can find that the intersection of the hypersurfaces $H_i^{(R)}=0$ contains a linear space $L^{(R)}$ of dimension $\lambda^{(R)}$ as long as the number of variables exceeds
\begin{align*}
    \gamma_3\big( C_2(R) m^8, \lambda^{(R)}\big) \ll_R m^{48} + m^{24} (\lambda^{(R)})^3.
\end{align*}
Thus we can reduce the problem to solving the remaining $R-1$ equations on $L^{(R)}$. In order for the residual system to be accessible to our methods, we need $\lambda^{(R)} \geq \gamma_5(R-1, m)$, whence by the induction hypothesis $\gamma_5(R,m)$ is bounded by
\begin{align*}
    \gamma_5(R,m) &\ll_R  m^{48} + m^{24} \big(\gamma_5(R-1, m)\big)^3 \\
    & \ll_R  m^{48} + m^{24} \Big( m^{12(3^{R-2}-1)+48 \cdot 3^{R-2}} \Big)^3.
\end{align*}
We may therefore conclude that there exists some function $A(R)$ such that
    \[\gamma_5(R,m)  \leq A(R) m^{{12(3^{R-1}-1)+48 \cdot 3^{R-1}}}, \]
as claimed. In fact, one may verify that the statement holds with
    \[A(R) \leq (18R)^{43 \cdot 3^R}.\]

It should be mentioned that the methods used here do not yield honest linear spaces unless we can ensure that the span of the vectors is of the right dimension. Let $s>\gamma_d(R,m)$ for some $m$ and suppose for convenience here that $m$ is even. Then we can find distinct vectors $\B x_1, \ldots, \B x_{m}$ contained in the intersection of the $R$ hyperplanes, and standard arguments (see the proof of Theorem 3 in \cite{lsp}) show that the vectors $\B x_1, \ldots, \B x_{m}$ are at least of rank $m/2$, provided that $s > 2Rrd/m+m/2$. This requirement is, however, easily met by any $s$ obeying the statements.

It is natural to ask what one should expect to be the true lower bound on the number $s$ of variables that ensures the existence of linear spaces. A comparison with related problems in $\mathbb R$ and $\mathbb C$ yields bounds which seem to intimate that the true growth rate of $s$ might be proportional to $Rm^{d-1}$. (For the former see Theorem 4 of \cite{lsp} and Theorem 6 of \cite{rd08}. In the latter case, apart from results concerning the multilinear Waring's problem such as \cite{alg-geom-1}, there is a bound for the existence of linear spaces due to Langer \cite{langer}, which is sharp but considers only the case $R=1$.) This evidence is further corroborated by the shape of the exponent $ms-Rrd$ in the main term in Theorems~\ref{Hasse} and \ref{asymp}, which will be positive only if $s \gg Rm^{d-1}$.
On the other hand, upper bounds have been provided by Dietmann and Wooley \cite{rd-tdw}, who proved that in the case $d=3$ and $R=1$ a growth rate of  $s \gg m^2$ suffices to guarantee the existence of $m$-dimensional linear spaces. While these results strongly suggest that the true growth rate of $s$ in $m$ really is $\asymp m^{d-1}$, there is less evidence for the growth rate in $R$ and there might be anomalies that have as yet not been spotted. The condition on $s$ we have established in Theorem \ref{asymp} and consequently Theorems \ref{cubic} and \ref{quintic} obviously fall short of the expected values, mainly because of the relatively large contribution arising from the local solubility condition. However, even in Theorem~\ref{Hasse} we miss the aim by a factor of $mR$.

\end{document}